\documentclass{arxiv}

\usepackage[shortlabels]{enumitem}
\setlist{leftmargin=*}
\usepackage{amssymb,mathtools}
\usepackage[hidelinks]{hyperref}
\usepackage{cleveref}
\usepackage{siunitx}
\usepackage{tikz-cd}
\usepackage{lmodern}

\numberwithin{equation}{section}

\Crefname{thrm}{Theorem}{Theorems}
\Crefname{lmm}{Lemma}{Lemmas}
\Crefname{xmpl}{Example}{Examples}
\Crefname{crllr}{Corollary}{Corollaries}
\Crefname{dfntn}{Definition}{Definitions}
\Crefname{rmrk}{Remark}{Remarks}

\newcommand{\accone}{\hat}
\newcommand{\acctwo}{\breve}
\newcommand{\accthr}{\mathring}

\newcommand{\Vb}{\accone{V}}
\newcommand{\vb}{\accone{v}}
\newcommand{\nb}{\accone{n}}
\newcommand{\cb}{\accone{c}}
\newcommand{\db}{\accone{d}}

\newcommand{\ub}{\accone{u}}
\newcommand{\wb}{\accone{w}}
\newcommand{\pib}{\accone{\pi}}

\newcommand{\Vo}{\acctwo{V}}
\newcommand{\uo}{\acctwo{u}}

\newcommand{\no}{\acctwo{n}}
\newcommand{\pio}{\acctwo{\pi}}

\newcommand{\V}{V}
\newcommand{\Vz}{\accthr{V}}
\newcommand{\Ph}{\mathfrak{H}}

\newcommand{\uz}{\accthr{u}}
\newcommand{\nz}{\accthr{n}}
\newcommand{\rhoz}{\rho_0}
\newcommand{\gamz}{\accthr{\gamma}}
\newcommand{\piz}{\accthr{\pi}}

\newcommand{\Pb}[1][1]{\accone{P}_{#1}^- \Lambda}
\newcommand{\Po}[1][1]{\acctwo{P}_{#1}^- \Lambda}
\newcommand{\Pz}[1][1]{\accthr{P}_{#1}^- \Lambda}
\newcommand{\PL}[1][1]{P_{#1}^- \Lambda}

\DeclareMathOperator{\Ker}{Ker}
\DeclareMathOperator{\Ran}{Ran}

\DeclareMathOperator{\grad}{grad}
\DeclareMathOperator{\curl}{curl}
\DeclareMathOperator{\divi}{div}
\DeclareMathOperator{\rot}{rot}

\begin{document}

\title{Solvers for mixed finite element problems using Poincaré operators based on spanning trees}

\author{Wietse M. Boon}
\address{University of Duisburg-Essen, Thea-Leymann Str. 9, Germany, \email{wietse.boon@uni-due.de}}



\begin{abstract}
	We propose a decomposition of Hilbert complexes that directly leads to a Poincaré operator. An explicit example is provided that decomposes a finite element differential complex using spanning trees in the grid. The Poincaré operator has three implications. First, it yields a new basis in which the mixed formulation of the Hodge-Laplace problem unravels from a large saddle point system into seven smaller, symmetric positive definite systems. These systems can be solved sequentially, and three of these have the same dimensions as the cohomology classes.
	Second, we use the operator to construct an explicit basis for the harmonic forms. Third, we propose an auxiliary space preconditioner for problems in weighted Sobolev spaces, that robustly handles the large kernel of the differential operator. These three implications are validated through numerical experiments.
\end{abstract}

\maketitle

\section{Introduction}
\label{sec: Introduction}

Mixed finite element methods are capable of preserving the structure of partial differential equations and thereby mimic physical behavior of a system at the discrete level. Stability and convergence of these methods can be derived rigorously using underlying theory \cite{arnold2018finite,boffi2013mixed}. However, the incorporation of physical conservation laws often leads to large saddle point problems that are computationally demanding to solve. Thus, efficient numerical solvers are essential to keep simulation times manageable. In this work, we propose such solvers by leveraging spanning trees in the grid.

The use of spanning trees in mixed finite element discretizations and solvers dates back, at least, to the tree-cotree decomposition of the Raviart-Thomas space in \cite{alotto1999mixed}. Later, they were employed in \cite{hiptmair2002generators,rodriguez2013construction,alonso2018efficient} to construct representatives of the first and second cohomology groups of the de Rham complex. In \cite{alonso2017graphs}, a spanning tree was employed to construct divergence-free finite element spaces. The key idea there is to restrict an $H(\curl)$-conforming finite element space by setting the edge degrees of freedom on a spanning tree to zero. This creates a subspace that is linearly independent of the kernel of the $\curl$ and, in turn, its image under the $\curl$ forms a basis for a divergence-free subspace of $H(\divi)$. This strategy was then generalized to higher-order spaces in \cite{alonso2018graph,devloo2022efficient,rodriguez2024basis}.
We note that in \cite{alonso2018graph}, the authors moreover employed a spanning tree of the dual graph to compute finite element functions with assigned divergence. This same strategy was later employed in \cite{boon2023reduced,boon2025deep} to ensure mass conservation in reduced order models.

In this paper, we build on these results by casting spanning tree decompositions in the general framework of Hilbert complexes. We identify sufficient conditions for such a decomposition to permit an explicit Poincaré operator.
These Poincaré operators are then used to obtain a new decomposition into subspaces related to the image and domain of the differential operators, and the harmonic forms. This provides a novel perspective on the constructions from \cite{alonso2017graphs,alonso2018graph,devloo2022efficient} as special cases of permitting decompositions in the context of Finite Element Exterior Calculus \cite{arnold2006finite,arnold2018finite} and the finite element de Rham complex.

We highlight three practical implications that arise from the derived theory. First, the decomposition provides a different basis for the Hilbert complex. Expressed in this basis, the Hodge-Laplace problem becomes a series of seven smaller, symmetric positive definite systems that can be solved sequentially. We show numerically that this results in a significant speed-up for all Hodge-Laplace problems on the complex, especially for the problem posed on the edge-based Nédélec space. The same solution is recovered, so this procedure forms a direct solver.

Second, the Poincaré operator allows us to construct an explicit basis for the harmonic forms. Each basis function requires solving two linear systems, and the numerically obtained basis functions are in qualitative agreement with expectation.

Third, we fit the new basis in the framework of auxiliary space preconditioning to construct a preconditioner for a more general class of problems. We show that the performance of this preconditioner is related to a Poincaré constant on a subspace of the decomposition. This Poincaré constant is computable and forms an upper bound on the Poincaré constant of the full finite element space. Numerical experiments indicate that the preconditioner is particularly well-suited for problems in which the differential operators dominate, which implies that it properly handles the large kernel of the differential operator.

We mention three further ties of this work to the existing literature. The decomposition is reminiscent of the regular decomposition from \cite{hiptmair2007nodal}. Additionally, this work draws inspiration from \cite{vcap2023bounded}, in which Poincaré operators are constructed for the twisted and BGG complexes. Finally, we refer the interested reader to \cite{ern2025discrete} for a recent review on discrete Poincaré operators.

The remainder of this paper is organized as follows. First, \Cref{sub: preliminaries} introduces preliminary concepts such as Hilbert complexes and the Whitney forms. \Cref{sec: permitting decomposition} introduces our main tool, mainly a \emph{permitting decomposition} of the Hilbert complex, along with several examples. The main theoretical results are presented in \Cref{sec: the Poincare operator}, where we construct a Poincaré operator that leads to a new decomposition of the Hilbert spaces. Using this decomposition, we propose a solver for the Hodge-Laplace problem in \Cref{sec: HL solver}. In \Cref{sec: aux preconditioning}, we then construct a preconditioner for finite element problems in weighted Sobolev spaces by placing the new decomposition in the framework of auxiliary space preconditioning. We present numerical experiments in \Cref{sec: Numerical results} that highlight the potential of these approaches. \Cref{sec: Concluding remarks} contains concluding remarks.

\subsection{Preliminary definitions and notation conventions}
\label{sub: preliminaries}

Let $(W^\bullet, d^\bullet)$ be a Hilbert complex consisting of a sequence of Hilbert spaces $W^k$, indexed by $k \in \mathbb{Z}$, and densely defined unbounded linear mappings $d^k : W^k \to W^{k + 1}$, with the property $\Ran(d^k) \subseteq \Ker(d^{k + 1})$, i.e.~$d^{k + 1}d^k = 0$ for all $k$. $(\cdot, \cdot)_W$ denotes the $W^k$-inner product, where we omit the superscript $k$ for brevity, and the associated norm is denoted as $\| \cdot \|$, in another slight abuse of notation.

Let $(V^\bullet, d^\bullet)$ be a subcomplex on which each $d^k$ is bounded. In this work, $V^k$ will either be the domain of $d^k$ or a finite-dimensional subspace, i.e. a conforming finite element space. These spaces are endowed with the graph norm
\begin{align} \label{eq: V norm}
	\| u \|_{V^k}^2 \coloneqq \| u \|^2 + \| d^k u \|^2.
\end{align}

Since $(V^\bullet, d^\bullet)$ is itself a Hilbert complex, it can be illustrated as
\begin{equation} \label{eq: Hilbert complex}
	\begin{tikzcd}
		\cdots
		\arrow[r, "d^{k - 2}"]
		&
		V^{k - 1}
		\arrow[r, "d^{k - 1}"]
		&
		V^k
		\arrow[r, "d^k"]
		&
		V^{k + 1}
		\arrow[r, "d^{k + 1}"]
		&
		\cdots
	\end{tikzcd}
\end{equation}
By the defining properties, two consecutive steps in this diagram map to zero.
The complex is called \emph{exact} if $\Ran(d^k) \supseteq \Ker(d^{k + 1})$ holds as well, i.e. if the range and kernel spaces coincide. In general, we define the space of harmonic $k$-forms as
\begin{align}
	\Ph^k = \{ u \in V^k : d^k u = 0 \text{ and } (u, d^{k - 1} v')_W = 0,\ \forall v' \in V^{k - 1} \},
\end{align}
which are trivial if the complex is exact. For simplicity, we assume that $(V^\bullet, d^\bullet)$ has the compactness property (cf.~\cite[Def.~4.3]{arnold2018finite}) and thus each $\Ph^k$ is finite-dimensional.

\begin{xmpl}
	Our main example concerns the \emph{de Rham complex}\cite[Sec.~4.3]{arnold2018finite}, which we describe next. For $n \in \{2, 3\}$, let $\Omega \subset \mathbb{R}^n$ be a Lipschitz domain. $L^2(\Omega)$ denotes the space of square integrable functions on $\Omega$. Let $(u, u')_\Omega$ be the inner product between $u, u' \in L^2(\Omega)$ and $\| u \|_\Omega^2 \coloneqq (u, u)_\Omega$ the corresponding norm. The base spaces are then given by $W^0 = W^n = L^2(\Omega)$ and $W^k = L^2(\Omega) \otimes \mathbb{R}^n$ for $0 < k < n$.

	In 3D, the differential operators correspond to $d^0 \coloneqq \grad$, $d^1 \coloneqq \curl$, and $d^2 \coloneqq \divi$. In 2D, we have $d^0 \coloneqq \rot = (\partial_y, -\partial_x)$ and $d^1 \coloneqq \divi$. For both cases of $n$, we let $d^n \coloneqq 0$. For $0 \le k \le n$, we set $V^k := H \Lambda^k$ as the subspace of square integrable, (and vector-valued for $0 < k < n$,) functions on $\Omega$ that are finite in the norm
	\begin{align} \label{eq: H norm}
		\| u \|_{H \Lambda^k}^2 \coloneqq \| u \|_\Omega^2 + \| d^k u \|_\Omega^2.
	\end{align}
	For completeness, we define $d^k \coloneqq 0$ and $W^k \coloneqq \{0\}$ for $k \not \in [0, n]$.
\end{xmpl}

We continue by introducing a key concept for this work, namely Poincaré operators.
\begin{dfntn} \label[dfntn]{def: poincare operator}
	A \emph{Poincaré operator} for a Hilbert complex $(V^\bullet, d^\bullet)$ is a family of linear maps $p_k : V^k \to V^{k - 1}$ that satisfies the relation
	\begin{align} \label{eq: poincare identity intro}
		d^{k - 1} p_k + p_{k + 1} d^k = I - \rhoz^k,
	\end{align}
	for all $k$, with $\rhoz^k: \V^k \to \Ker(d^k)$ a linear map.
\end{dfntn}

This work focuses on two model problems. First, the mixed variational formulation of the Hodge-Laplace problem on $V^k$ reads: given functionals $g \in (V^{k - 1})^*$ and $f \in (V^k)^*$, find $(v, u) \in V^{k - 1} \times V^k$ such that
\begin{subequations}  \label{eq: HL problem cont}
	\begin{align}
		(v, v')_W - (u, d^{k - 1} v')_W                        & = \langle g, v' \rangle,    &
		\forall v'                                             & \in V^{k - 1},                \\
		(d^{k - 1} v, u')_W + (d^k u, d^k u')_W - (\phi, u')_W & = \langle f, u' \rangle,    &
		\forall u'                                             & \in V^k                       \\
		(u, \phi')_W                                           & = \langle e, \phi' \rangle, &
		\forall \phi'                                          & \in \Ph^k.
	\end{align}
\end{subequations}
Here, and throughout this paper, we adopt the notation $V^*$ for the dual of a function space $V$ and $\langle \cdot, \cdot \rangle$ denotes a duality pairing.
Apostrophes denote test functions.

Second, we consider the projection problem onto $V^k$, which reads: for given $f \in (V^k)^*$ and $\alpha > 0$, find $u \in V^k$ such that
\begin{align} \label{eq: projection problem cont}
	\alpha^2 (u, u')_W + (d^k u, d^k u')_W
	           & = \langle f, u' \rangle, &
	\forall u' & \in V^k.
\end{align}

\section{Permitting decompositions of Hilbert complexes}
\label{sec: permitting decomposition}

In this section, we introduce the theoretical tools that serve as the building blocks for the proposed numerical solvers. We first introduce a particular decomposition of the spaces $\V^k$ in \Cref{sub: permitting decomposition}. \Cref{sub: Spanning trees} concerns our main example, namely a spanning tree decomposition of the Whitney forms. This construction is then generalized to several additional examples in \Cref{sub: additional examples}.

\subsection{A permitting decomposition}
\label{sub: permitting decomposition}

We consider decompositions of the form $\V^\bullet = \V_1^\bullet \oplus \V_2^\bullet$, which means that for each $k$, each $v \in V^k$ decomposes uniquely as $v = v_1 + v_2$ with $(v_1, v_2) \in V_1^k \times V_2^k$.

Let $\Vb^k$, $\Vo^k$, and $\Vz^k$ be subspaces of $\V^k$, endowed with projection operators $\pib_k: \V^k \to \Vb^k$, $\pio_k: \V^k \to \Vo^k$, and $\piz_k: \V^k \to \Vz^k$. We assume that the projection operators vanish outside the corresponding subspace, i.e.~$\pib_k (\Vo^k \oplus \Vz^k) = 0$ etc.
We omit subscripts on the projections for brevity, as these will be clear from context.
Moreover, we define $\db^k: \Vb^k \to \V^{k + 1}$ as the restriction of $d^k$ on $\Vb^k$.
With these definitions in place, we introduce the concept of a permitting decomposition.

\begin{dfntn} \label[dfntn]{def: permitting}
	A decomposition of $(\V^\bullet, d^\bullet)$ given by $\V^\bullet = \Vb^\bullet \oplus \Vo^\bullet \oplus \Vz^\bullet$ is \emph{permitting} if the following two properties are satisfied:
	\begin{enumerate}
		\item
		      The operator
		      \begin{align}
			      \pio \db^k : \Vb^k \to \Vo^{k + 1}
		      \end{align}
		      is invertible for all $k$.
		\item
		      For each $k$, a projection $\gamz: \V^k \to \Vz^k$ exists that is a bijection between $\Ph^k$ and $\Vz^k$ and satisfies
		      \begin{align} \label{eq: permitting 2}
			      \Vb^k           & \subseteq \Ker(\gamz), &
			      \gamz d^{k - 1} & =0.
		      \end{align}
	\end{enumerate}
\end{dfntn}

A permitting decomposition of $(\V^\bullet, d^\bullet)$ can be illustrated as
\begin{equation}
	\begin{tikzcd}[column sep=large]
		\cdots
		\arrow[dr, dashed]
		&
		\Vb^{k - 1}
		\arrow[d, "\oplus" description, line width = 0mm]
		\arrow[dr, "\pio \db^{k - 1}" description]
		&
		\Vb^k
		\arrow[d, "\oplus" description, line width = 0mm]
		\arrow[dr, "\pio \db^k" description]
		&
		\Vb^{k + 1}
		\arrow[d, "\oplus" description, line width = 0mm]
		\arrow[dr, dashed]
		\\[-10pt]
		&
		\Vo^{k - 1}
		\arrow[d, "\oplus" description, line width = 0mm]
		&
		\Vo^k
		\arrow[d, "\oplus" description, line width = 0mm]
		&
		\Vo^{k + 1}
		\arrow[d, "\oplus" description, line width = 0mm]
		&
		\cdots
		\\[-10pt]
		&
		\Vz^{k - 1}
		&
		\Vz^k
		&
		\Vz^{k + 1}
		&
	\end{tikzcd}
\end{equation}
in which the operators on the diagonal lines are invertible and the spaces $\Vz^k$ are isomorphic to $\Ph^k$, due to the operator $\gamz$.

\subsection{Main example: the Whitney forms}
\label{sub: Spanning trees}

In this section, we introduce a family of low-order finite element spaces, known as the Whitney forms, that form a discrete Hilbert complex. We then construct a decomposition based on spanning trees, borrowing concepts from \cite{alonso2017graphs,devloo2022efficient,hiptmair2002generators}. Its permitting property is shown afterward in \Cref{lem: permitting}.

\begin{dfntn}[Whitney forms] \label[dfntn]{def: whitney}
	Let $\Omega_h$ be a simplicial grid that tesselates a connected Lipschitz domain $\Omega \subset \mathbb{R}^n$ with $n \in \{2, 3\}$. We denote the Betti numbers of $\Omega$ by $b_k$, which correspond to the numbers of $k$-holes in the domain. Let $\PL^k$ be the trimmed finite element space of lowest order\cite[Sec.~7.6]{arnold2018finite}, with one degree of freedom per $k$-simplex of the grid $\Omega_h$. These spaces are given by
	\begin{itemize}
		\item $\PL^0 \coloneqq \mathbb{L}_1(\Omega_h)$: the nodal-based, linear Lagrange finite element space.

		\item If $n = 3$, then $\PL^1 \coloneqq \mathbb{N}_0(\Omega_h)$: the lowest-order, edge-based Nédélec space of the first kind \cite{nedelec1980mixed}.

		\item $\PL^{n - 1} \coloneqq \mathbb{RT}_0(\Omega_h)$: the facet-based Raviart-Thomas space of lowest order \cite{raviart-thomas-0}.

		\item $\PL^n \coloneqq \mathbb{P}_0(\Omega_h)$: the space of element-wise constants.
	\end{itemize}
\end{dfntn}

\begin{figure}[ht]
	\centering
	\includegraphics[width = 0.40\textwidth]{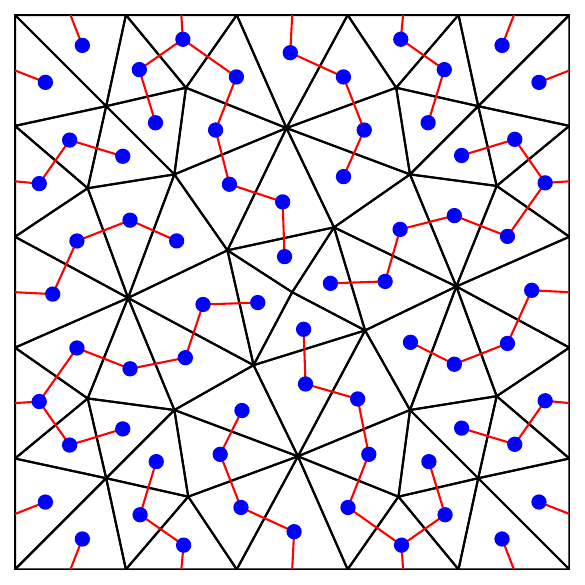}
	\includegraphics[width = 0.40\textwidth]{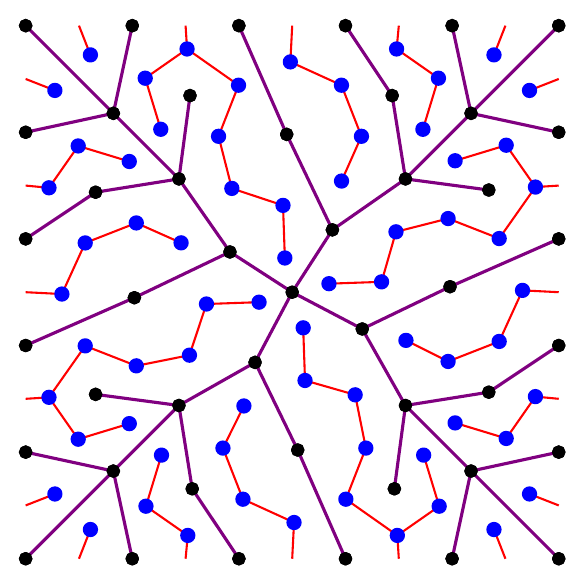}
	\caption{Decomposition of the facet-based space $\PL^1$ on an unstructured grid of a topologically trivial 2D domain. (Left) A spanning tree on the dual graph connects all cells, including a fictitious ``outside'' cell through the domain boundary. The subspace $\Pb^1$ has degrees of freedom on the facets crossed by this tree. (Right) The complementary set of facets $\acctwo{\mathcal{F}}$, illustrated in purple, forms a spanning tree for the nodes. $\Po^1$ is the subspace of $\PL^1$ associated with these facets.}
	\label{fig: trees}
\end{figure}

\begin{figure}[ht]
	\centering
	\includegraphics[width = 0.40\textwidth]{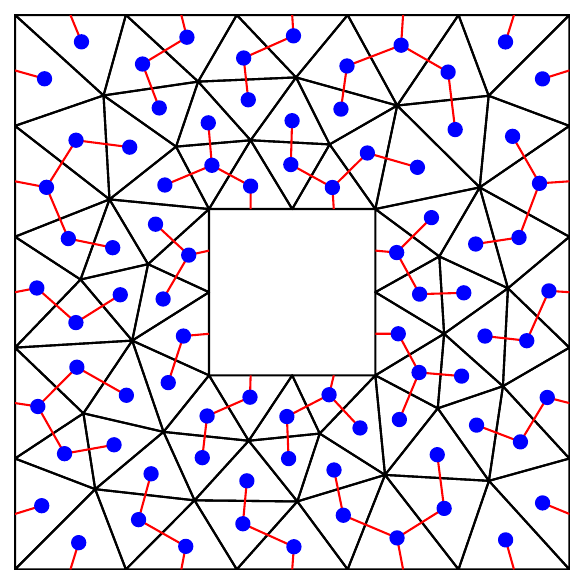}
	\includegraphics[width = 0.40\textwidth]{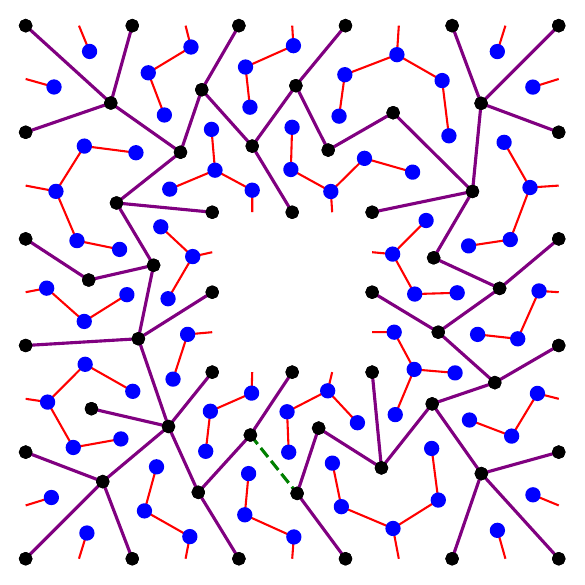}
	\caption{The decomposition of $\PL^1$ in the case of non-trivial topology. A third subspace $\Pz^1$ is introduced that has a degree of freedom on a single mesh facet. This facet, indicated with a dashed green line, is located in the central bottom of the domain and is the only member of the set $\accthr{\mathcal{F}}$. A closed curve around the hole is contained in $\acctwo{\mathcal{F}} \cup \accthr{\mathcal{F}}$.}
	\label{fig: holy trees}
\end{figure}

\begin{xmpl}[A spanning tree decomposition] \label[xmpl]{ex: spanning tree decomposition}
	We construct the decomposition of $(\PL^\bullet, d^\bullet)$ for each of the four non-trivial cases, in descending order
	\begin{itemize}
		\item $k = n$. For the $n$-forms, we set
		      \begin{align} \label{eq: n-form decomp}
			      \Pb^n & \coloneqq \{0\}, & \Po^n & \coloneqq \PL^n, & \Pz^n & \coloneqq \{0\}, & \gamz & = 0.
		      \end{align}
		\item $k = n - 1$.
		      Let $\mathcal{G}_{n - 1}^*$ be the dual graph whose vertices correspond to the cells of $\Omega_h$ with an additional vertex $m^\infty$ associated with the ``outside'' $\mathbb{R}^n \setminus \Omega_h$. An edge on $\mathcal{G}_{n - 1}^*$ implies that the two corresponding cells share a face. Thus, all boundary faces on $\partial \Omega_h$ are represented on $\mathcal{G}_{n - 1}^*$ by edges adjacent to $m^\infty$. Let $\mathcal{T}_{n - 1}^*$ be a spanning tree of $\mathcal{G}_{n - 1}^*$ and let $\accone{\mathcal{F}}$ be the set of faces in $\Omega_h$ that correspond to the edges of $\mathcal{T}_{n - 1}^*$.

		      We initialize $\acctwo{\mathcal{F}}$ as the set complementary to $\accone{\mathcal{F}}$ and $\accthr{\mathcal{F}} = \emptyset$. For the next construction, we introduce $\mathcal{S}$ as a copy of $\acctwo{\mathcal{F}}$. We then prune $\mathcal{S}$ by iteratively removing all faces in the set that are connected to less than $n$ neighbors in $\mathcal{S}$. If $\Omega$ is topologically trivial, then this process annihilates the set, cf.~\Cref{fig: trees}. In 2D, this is a consequence of the tree-cotree decomposition of graphs \cite[Thm.~XI.6.]{tutte2001graph}.

		      Otherwise, the pruning process will leave a closed surface (curve in 2D) that surrounds the $(n - 1)$-holes of the domain, cf.~\Cref{fig: holy trees}. There are $b_{n - 1}$ of such holes, corresponding to the $(n - 1)$-th Betti number. The remaining surface divides the domain into $b_{n - 1} + 1$ subdomains, which we denote by $\Omega_i$. We again create a dual graph, this time letting the vertices correspond to the subdomains and the edges indicating neighboring subdomains. The edges of the spanning tree on this graph correspond to $b_{n - 1}$ subdomain interfaces. On each of the interfaces corresponding to that tree, we select a single facet and move it from $\acctwo{\mathcal{F}}$ to $\accthr{\mathcal{F}}$.

		      Using these facet sets, the subspaces of $\PL^{n - 1}$ are defined as
		      \begin{align*}
			      \Pb^{n - 1} & \coloneqq \PL^{n - 1}(\accone{\mathcal{F}}), &
			      \Po^{n - 1} & \coloneqq \PL^{n - 1}(\acctwo{\mathcal{F}}), &
			      \Pz^{n - 1} & \coloneqq \PL^{n - 1}(\accthr{\mathcal{F}}).
		      \end{align*}
		      where we employ the short-hand notation
		      \begin{align}
			      \PL^{n - 1}(\mathcal{F})\ \coloneqq \left\{ u \in \PL^{n - 1} : \int_f \nu_f \cdot u = 0, \ \forall f \not \in \mathcal{F} \right\},
		      \end{align}

		      with $\nu_f$ the normal vector on facet $f$.
		      The operator $\gamz$ is then defined such that
		      \begin{align}
			      \gamz                                              & : \PL^{n - 1} \to \Pz^{n - 1}, &
			      \int_{\partial \Omega_i} \nu_i \cdot (\gamz u - u) & = 0,                           &
			      \forall \Omega_i,
		      \end{align}
		      with $\nu_i$ the unit normal vector on $\partial \Omega_i$.
		\item $k = 1$, $n = 3$.
		      In this case, we first identify the closed loops in the grid that surround $1$-holes using the strategy proposed in \cite{hiptmair2002generators} on the boundary $\partial \Omega_h$. The construction involves creating a tree-cotree decomposition on the boundary to identify a class of loops. After appropriate shifting these curves, we compute linking numbers to identify the relevant loops that surround $1$-holes, which we denote by $\Gamma_i$. We refer to \cite{hiptmair2002generators} for further details. We then choose one edge per loop and use these to form the set $\accthr{\mathcal{E}}$.

		      Next, we interpret the vertex and edge sets of $\Omega_h$ as a graph $\mathcal{G}_1$. $\mathcal{T}_1$ is then constructed as a spanning tree of $\mathcal{G}_1$, under the constraint that the edges on $\Gamma_i \setminus \accthr{\mathcal{E}}$ belong to $\mathcal{T}_1$ for each loop $\Gamma_i$ identified above. $\acctwo{\mathcal{E}}_1$ is then given by the edges of $\mathcal{T}_1$.

		      With the complementary set $\acctwo{\mathcal{E}} = \Delta_1 \setminus (\acctwo{\mathcal{E}} \cup \accthr{\mathcal{E}})$, we define the subspaces as
		      \begin{align*}
			      \Pb^1 & \coloneqq \PL^1 (\accone{\mathcal{E}}), &
			      \Po^1 & \coloneqq \PL^1 (\acctwo{\mathcal{E}}), &
			      \Pz^1 & \coloneqq \PL^1 (\accthr{\mathcal{E}}).
		      \end{align*}
		      in which
		      \begin{align*}
			      \PL^1(\mathcal{E}) & \coloneqq \left\{ u \in \PL^1 : \int_e t_e \cdot u = 0, \ \forall e \not \in \mathcal{E} \right\},
		      \end{align*}
		      with $t_e$ the tangent vector on edge $e$. In this case, the operator $\gamz$ is such that
		      \begin{align}
			      \gamz                                   & : \PL^1 \to \Pz^1, &
			      \int_{\Gamma_i} t_i \cdot (\gamz u - u) & = 0,               &
			      \forall \Gamma_i,
		      \end{align}
		      with $t_i$ the tangent vector to $\Gamma_i$.

		\item $k = 0$. Let $m^0$ be an arbitrarily chosen node from the node set $\Delta_0$. We decompose the space into the functions that are zero on $m^0$ and the functions that are zero on the other nodes:
		      \begin{align*}
			      \Pb^0 & \coloneqq \PL^0(\Delta_0 \setminus m^0), &
			      \Po^0 & \coloneqq \{0\},                         &
			      \Pz^0 & \coloneqq \PL^0(m^0).
		      \end{align*}
		      The operator $\gamz$ is simply defined as
		      \begin{align}
			      \gamz        & : \PL^0 \to \Pz^0, &
			      \gamz u(m^0) & = u(m^0).
		      \end{align}
	\end{itemize}

	Finally, the projections $\pib$, $\pio$, and $\piz$ are defined as restriction operators onto the corresponding subspaces.
\end{xmpl}

Possible realizations of \Cref{ex: spanning tree decomposition} on two-dimensional grids are illustrated in \Cref{fig: trees,fig: holy trees}. Here, we have generated $\mathcal{T}_1^*$ by using a breadth-first search algorithm on the dual graph, rooted at the outside vertex $m^\infty$. From these figures, we observe that the space $\PL^1$ can be decomposed in degrees of freedom associated with a map to the cells (the divergence), a mapping from the nodes (rot), and a subspace with the same dimension as the space of harmonic forms. This observation is made more precise in the remainder of this article.

\begin{lmm} \label[lmm]{lem: permitting}
	The decomposition from \Cref{ex: spanning tree decomposition} is permitting in the sense of \Cref{def: permitting}.
\end{lmm}
\begin{proof}
	We start with the second condition of \Cref{def: permitting} and consider $\gamz$ in the three non-trivial cases.
	\begin{itemize}
		\item $k = n - 1$. Each $\gamz$ measures the normal component of the vector field through the surfaces that surround $(n - 1)$-holes. The harmonic forms in this case can be interpreted as divergence-free flux fields that transport mass from one closed boundary of $\Omega$ to another. Thus, if $\gamz \phi = 0$ for $\phi \in \Ph$, then $\phi = 0$. Moreover, $\gamz d^{n - 2} v = 0$ for any $v \in \V^{n - 2}$ because the integral over a closed surface of $d^{n - 2} v$ is zero by Stokes theorem. Finally, the decomposition is chosen such that $\Vb^{n - 1}$ does not have degrees of freedom on these surfaces, so $\Vb^{n - 1} \subseteq \Ker(\gamz)$.
		\item $k = 1$, $n = 3$. In this case, $\gamz$ integrates the tangential component of a vector field along a closed curve that surrounds a $1$-hole. The harmonic forms in this case are exactly the circulations around these holes, and thus $\gamz \phi$ is only zero for $\phi \in \Ph^1$ if $\phi = 0$. Again, Stokes theorem and the chosen decomposition give us \eqref{eq: permitting 2}.
		\item $k = 0$. The space of harmonic forms is given by the constants ($\Ph^0 = \mathbb{R}$) and $\gamz u$ is invertible on this space. \eqref{eq: permitting 2} follows immediately from the definitions.
	\end{itemize}

	For the first condition, we need to show that $\pio \db^k : \Pb^k \to \Po^{k + 1}$ is invertible for each $k$. Since $\Pb^n$ is zero in this example, we are left with three non-trivial cases to consider:
	\begin{itemize}
		\item $k = 0$. In this case, the operator $\pio \db$ maps from the nodes to the edges of the spanning tree $\mathcal{T}_1$. Let $|\Delta_0|$ be the number of nodes in $\Omega_h$, then $\dim( \Pb^0) = |\Delta_0| - 1$ because of the removed node $m^0$. The number of edges on the tree is also $|\Delta_0| - 1$ and so $\dim( \Pb^0) = \dim (\Po^1)$.

		      The operator $\pio \db^0$ can thus be represented as a square matrix and it remains to show that it is injective. Let $\ub \in \Pb^0$ satisfy $\pio \db^0 \ub = 0$. This implies that the difference between each node connected by the tree $\mathcal{T}_1$ is zero. Since $\mathcal{T}_1$ is a spanning tree, this means that $\ub$ is constant. The only constant in $\Pb^0$ is zero, so $\pio \db$ is injective.

		\item $k = n - 1$. Again, we first look at the dimensions of $\Pb^{n - 1}$ and $\Po^n$. Let $|\Delta_{n - 1}|$ and $|\Delta_n|$ denote the number of facets and cells in $\Omega_h$, respectively. Recall that the facet set $\mathcal{F}_{n - 1}$ corresponds to the edges of a spanning tree with $|\Delta_n| + 1$ vertices, so $\dim( \Pb^{n - 1}) = |\Delta_n|$. Since $\Po^n = \PL^n$, we also have $\dim( \PL^n) = |\Delta_n|$.

		      It remains to show injectivity. Let $\ub \in \Pb^{n - 1}$ with $\pio \db^{n - 1} \ub = 0$. Since $\Po^n = \PL^n$, this implies that $d^{n - 1} \ub = 0$. Now consider the cells associated with the leaves of the spanning tree $\mathcal{T}_{n - 1}^*$. On each of these cells, $d^{n - 1} \ub = \divi \ub$ is determined by one degree of freedom of $\ub$, which must therefore be zero. We may then prune the tree by removing these leaves and repeat the argument on the reduced tree. Iterating the argument until we reach the root of the tree, we conclude $\ub = 0$.

		\item $k = 1$, $n = 3$. For the dimensions of $\Pb^1$ and $\Po^2$, let $|\Delta_1|$ denote the number of edges in $\Omega_h$. Using the cases $k = 0$ and $k = n - 1$ above, we derive
		      \begin{align}
			      \dim( \Pb^1) & = \dim( \PL^1) - \dim( \Po^1) - \dim( \Pz^1) = |\Delta_1| - (|\Delta_0| - 1) - b_1, \\
			      \dim( \Po^2) & = \dim( \PL^2) - \dim( \Pb^2) - \dim( \Pz^2) = |\Delta_2| - |\Delta_3| - b_2
		      \end{align}

		      Considering $\Omega_h$ as a pure simplicial complex, we use the following relation between the Euler characteristic and the Betti numbers
		      \begin{align}
			      \chi(\Omega_h) \coloneqq |\Delta_0| - |\Delta_1| + |\Delta_2| - |\Delta_3| = b_0 - b_1 + b_2 - b_3
		      \end{align}

		      Using the fact that $b_0 = 1$ and $b_3 = 0$, we derive
		      \begin{align*}
			      \dim( \Pb^1) - \dim( \Po^2)
			      = |\Delta_1| - |\Delta_0| + 1 - b_1 - |\Delta_2| + |\Delta_3| + b_2
			      = \chi(\Omega_h) - \chi(\Omega_h) = 0
		      \end{align*}

		      It remains to show that $\pio \db^1$ is injective, so let $\ub \in \Pb^1$ satisfy $\pio \db^1 \ub = 0$, and thus $d^1 \ub \in \Pb^2 \oplus \Pz^2$.
		      Now $\gamz$ is identity on $\Pz^2$ and, by the calculation for $k = n - 1$ above, we have that $\pio \db^2$ is invertible on $\Pb^2$. Using this in combination with $\gamz d^1 = 0$ and $d^2 d^1 = 0$, we obtain
		      \begin{align*}
			      d^1 \ub
			      = (\pib + \piz) d^1 \ub
			      = (\pio \db^2)^{-1} (\pio \db^2) \pib d^1 \ub + \gamz \piz d^1 \ub = 0.
		      \end{align*}

		      Now, the Hodge decomposition on $(\PL^\bullet, d^\bullet)$ implies that $\ub = d^0 v + \phi$ with $v \in \PL^0$ and $\phi \in \Ph^1$.
		      Applying the operator $\gamz$ to both sides, we note that $\ub \in \Ker(\gamz)$ and $\gamz d^0 = 0$ so
		      \begin{align}
			      \gamz \phi = 0
		      \end{align}
		      But $\gamz$ is a bijection on $\Ph^1$, so $\phi = 0$.
		      We now choose $\vb \in \Pb^0$ by subtracting a constant from $v$. Since $\pio \db^0$ is invertible on $\Pb^0$, as shown above, we obtain
		      \begin{align}
			      \vb = (\pio \db^0 )^{-1} \pio d^0 \vb = (\pio \db^0 )^{-1} \pio \ub = 0.
		      \end{align}
		      Since both $\phi$ and $d^0 v$ are zero, we conclude that $\ub = 0$.
	\end{itemize}
	In summary, $\pio \db$ is invertible for all $k$ and so the decomposition from \Cref{ex: spanning tree decomposition} is permitting.
\end{proof}

\subsection{Additional examples}
\label{sub: additional examples}

\begin{xmpl}
	Let $(\PL[r]^\bullet, d^\bullet)$ denote the general complex of trimmed finite elements, with \Cref{def: whitney} corresponding to $r = 1$. In 3D, this complex is given by
	\begin{equation} \label{eq: P1 spaces}
		\begin{tikzcd}
			\mathbb{L}_{r}
			\arrow[r, "\grad"]
			&
			\mathbb{N}_{r - 1}
			\arrow[r, "\curl"]
			&
			\mathbb{RT}_{r - 1}
			\arrow[r, "\divi"]
			&
			\mathbb{P}_{r - 1}
		\end{tikzcd}
	\end{equation}

	$(\PL^\bullet, d^\bullet)$ is a subcomplex of each $(\PL[r]^\bullet, d^\bullet)$, with isomorphic cohomology classes, so we set each $\Pz[r]^k \coloneqq \Pz^k$ and each $\gamz$ exactly as in \Cref{ex: spanning tree decomposition}. We once again present the decomposition in descending order.
	\begin{itemize}
		\item $k = n$. Since $d^n = 0$ and $b_n = 0$, we set
		      \begin{align}
			      \Pb[r]^n & \coloneqq \{0\}, & \Po[r]^n & \coloneqq \PL[r]^n.
		      \end{align}
		\item $k = n - 1$. Let $\PL[r](\Delta_n)$ and $\PL[r](\Delta_{n - 1})$ denote the basis functions of $\mathbb{RT}_r$ with cell-based and facet-based degrees of freedom, respectively. Since $d^{n - 1} \PL[r](\Delta_n) = \mathbb{P}_r / \mathbb{P}_0$, we include these in the $\Pb^k$ space:
		      \begin{align*}
			      \Pb[r]^{n - 1} & \coloneqq \Pb^{n - 1} \oplus \PL[r](\Delta_n),                                               &
			      \Po[r]^{n - 1} & \coloneqq \Po^{n - 1} \oplus (\PL[r]^{n - 1}(\Delta_{n - 1}) / \PL^{n - 1}(\Delta_{n - 1})).
		      \end{align*}
		      Here, the notation $P_r / P_0$ denotes the polynomials in $P_r$ that vanish on the degrees of freedom of $P_0$.
		\item $k = 1$, $n = 3$. We first define the space for $r = 2$, whose degrees of freedom are facet- and edge-based:
		      \begin{align*}
			      \Pb[2]^1 & \coloneqq \Pb^1 \oplus \PL[1]^1(\Delta_2),                     &
			      \Po[2]^1 & \coloneqq \Po^1 \oplus (\PL[1]^1(\Delta_1) / \PL^1(\Delta_1)).
		      \end{align*}
		      Note that the gradients of $\PL^0$ are captured by the edge degrees of freedom in $\Po[2]^1$, and the remainder forms the space $\Pb[2]^1$.
		      For higher values of $r$, this construction generalizes as
		      \begin{align*}
			      \Po[r]^1 & \coloneqq \Po^1
			      \oplus \bigoplus_{i = 1}^3 d^0 (\PL[r]^0(\Delta_i) / \PL^0(\Delta_i)), \\
			      \Pb[r]^1 & \coloneqq \PL[r]^1 / (\Po[r]^1 \oplus \Pz[r]^1).
		      \end{align*}
		\item $k = 0$. Let
		      \begin{align*}
			      \Pb[r]^0 & \coloneqq \Pb^0 \oplus (\PL[r]^0 / \PL^0), &
			      \Po[r]^0 & \coloneqq \{0\}.
		      \end{align*}
	\end{itemize}
	The permitting property of this decomposition can again be shown by comparing the dimensions of the subspaces and proving injectivity of $\pio \db^k$.
\end{xmpl}

\begin{xmpl}[Hodge decomposition] \label[xmpl]{ex: Hodge}
	We may alternatively define $\Po^k \coloneqq d^{k - 1} \PL^{k - 1}$, $\Pz^k \coloneqq \Ph^k$, and $\Pb^k$ as the $L^2$-orthogonal complement, i.e.~$\Pb^k \perp (\Po^k \oplus \Pz^k)$. This is the Hodge decomposition \cite[Sec.~4.2.2]{arnold2018finite}. The subspaces are then endowed with the $L^2$-projection operators $\pio$, $\pib$, and $\piz$. While a permitting decomposition, the construction of a basis for $\Pb^k$ and $\Po^k$ is not as straightforward as in \Cref{ex: spanning tree decomposition}. Moreover, the projection operators involve solving linear systems, which is more computationally demanding than applying restrictions. On a related note, the subspaces in \Cref{ex: spanning tree decomposition} are not orthogonal in the $W^k$-inner product, but the vector representations of their degrees of freedom are orthogonal in $\ell^2$.
\end{xmpl}

\begin{xmpl}[Virtual element methods] \label[xmpl]{ex: VEM}
	The construction from \Cref{ex: spanning tree decomposition} can be generalized to the case of mixed virtual element methods on polytopal grids \cite{beirao2013basic}. The lowest-order spaces $V_0 \Lambda^k$ have one degree of freedom for each $k$-dimensional polytope of the mesh and form an exact sequence $(V_0 \Lambda^\bullet, d^\bullet)$, as was shown in \cite{da2022virtual,boon2025nodal}. The arguments from \Cref{lem: permitting} transfer directly to this case with only one minor adjustment. The dual graph $\mathcal{G}_{n - 1}^*$ becomes a multi-graph because it is possible for a pair of elements to share multiple facets. Nevertheless, spanning trees can be constructed also for multi-graphs and the decomposition follows in the same manner.
\end{xmpl}

\section{The Poincaré operator}
\label{sec: the Poincare operator}

By the defining properties of a permitting decomposition, we are able to invert the operator $\pio \db^k$ for each $k$. This allows us to construct the following operators.

\begin{dfntn} \label[dfntn]{def: def p}
	Let $\V^\bullet = \Vb^\bullet \oplus \Vo^\bullet \oplus \Vz^\bullet$ be a permitting decomposition of $(\V^\bullet, d^\bullet)$. The \emph{permitted operators} $p_k : \V^k \to \V^{k -1}$ and $\rhoz^k : \V^k \to \V^k$ are given by 
	\begin{align}
		p_k     & \coloneqq (\pio \db^{k - 1})^{-1} \pio, &
		\rhoz^k & \coloneqq (I - p_{k + 1}d^k) \gamz,
	\end{align}
	for all $k$.
\end{dfntn}

We are now ready to present the main result of this work, in the following theorem.

\begin{thrm} \label{thm: poincare operator}
	The permitted operators from \Cref{def: def p} satisfy
	\begin{align} \label{eq: poincare identity}
		d^{k - 1}p_k + p_{k + 1}d^k = I - \rhoz^k,
	\end{align}
	for all $k$.
\end{thrm}
\begin{proof}
	Let $v \in \V^k$ and consider the decomposition $v = \gamz v + (I - \gamz) v$. For $\gamz v$, the definitions immediately give us $p_k \gamz = 0$ and so
	\begin{align} \label{eq: poincare z}
		(d^{k - 1}p_k + p_{k + 1}d^k) \gamz v
		= p_{k + 1}d^k \gamz v
		= (I - \rhoz^k) \gamz v.
	\end{align}

	It therefore suffices to consider $u = (I - \gamz)v$, which implies $\gamz u = 0$. We decompose $u = \ub + \uo + \uz$ such that $(\ub, \uo, \uz) \in \Vb^k \times \Vo^k \times \Vz^k$. We prove the equality by showing five identities.

	First, we note that for the first component $\ub$, we have $p_k \ub = 0$. This allows us to derive:
	\begin{align} \label{eq: ident 1}
		(d^{k - 1}p_k + p_{k + 1}d^k) \ub
		 & = ( (\pio \db^k)^{-1} \pio) d^k \ub
		= (\pio \db^k)^{-1} (\pio \db^k) \ub
		= \ub.
	\end{align}

	Second, since $\Ran(p_{k + 1}) \subseteq \Vb^k$, we have $\pio p_{k + 1} = 0$ and we obtain an identity concerning $\uo$:
	\begin{align} \label{eq: ident 2}
		\pio (d^{k - 1}p_k + p_{k + 1}d^k) \uo
		= \pio \db^{k - 1}p_k \uo
		 & = \pio \db^{k - 1} (\pio \db^{k - 1})^{-1} \pio \uo
		= \uo.
	\end{align}

	Third, the properties $p_k \uz = 0$, $\pio p_{k + 1} = 0 $, and $ \piz p_{k + 1} = 0$ provide
	\begin{align} \label{eq: ident 3}
		(\pio + \piz) (d^{k - 1}p_k + p_{k + 1}d^k) \uz & = 0.
	\end{align}

	The final two steps are slightly more involved. Let us recall that $0 = \gamz u = \gamz(\uo + \uz)$ and that $\gamz d^{k - 1} = 0$. This allows us to derive
	\begin{align*}
		\gamz \uz
		= \gamz (d^{k - 1} p_k \uo - \uo)
		 & = \gamz (\pio + \pib + \piz)(d^{k - 1} p_k \uo - \uo)                            \\
		 & = \gamz \pio(d^{k - 1} p_k \uo - \uo) + 0 + \gamz \piz (d^{k - 1} p_k \uo - \uo)
		= \gamz \piz d^{k - 1} p_k \uo,
	\end{align*}
	where the first term vanishes in the final equality due to \eqref{eq: ident 2}. Now note that both $\uz$ and $\piz d^{k - 1} p_k \uo$ are elements of $\Vz^k$. Since $\gamz$ is a bijection on this space, we conclude that
	\begin{align} \label{eq: uz is dp uo}
		\uz = \piz d^{k - 1} p_k \uo.
	\end{align}
	Together with $\piz p_{k + 1} = 0$, this provides the fourth identity
	\begin{align} \label{eq: ident 4}
		\piz (d^{k - 1}p_k + p_{k + 1}d^k) \uo
		= \uz.
	\end{align}

	For the fifth and final step, we observe that $p_{k + 1}$ is the composition of two surjective operators, and is therefore surjective onto $\Vb^k$, i.e.
	\begin{align} \label{eq: surjectivity}
		\Ran(p_{k + 1}) = \Vb^k.
	\end{align}

	In turn, a $\acctwo{t} \in \Vo^{k + 1}$ exists such that $p_{k + 1} \acctwo{t} = \pib d^{k - 1} p_k (\uo + \uz)$.
	We find $\acctwo{t}$ by applying the left-inverse operator $\pio \db^k : \Vb^k \to \Vo^{k + 1}$:
	\begin{align}
		\acctwo{t} 
		= \pio \db^k \pib d^{k - 1} p_k (\uo + \uz)
		 & = \pio d^k (I - \pio - \piz) d^{k - 1} p_k \uo \nonumber      \\
		 & = 0 - \pio d^k (\pio \db^{k - 1}) (\pio \db^{k - 1})^{-1} \uo
		- \pio d^k \piz d^{k - 1} p_k \uo \nonumber                      \\
		 & = -\pio d^k (\uo + \uz),
	\end{align}
	where the last equality is due to \eqref{eq: uz is dp uo}.
	With the derived $\acctwo{t}$, we obtain the fifth identity:
	\begin{align} \label{eq: ident 5}
		\pib (d^{k - 1}p_k + p_{k + 1}d^k) (\uo + \uz)
		 & = p_{k + 1} \acctwo{t} + p_{k + 1}d^k (\uo + \uz) \nonumber   \\
		 & = - p_{k + 1} \pio d^k (\uo + \uz) + p_{k + 1}d^k (\uo + \uz)
		= 0
	\end{align}

	The combination of the five identities \eqref{eq: ident 1}, \eqref{eq: ident 2}, \eqref{eq: ident 3}, \eqref{eq: ident 4}, and \eqref{eq: ident 5} now proves \eqref{eq: poincare identity}.
	%
\end{proof}

\Cref{thm: poincare operator} immediately implies several key properties.

\begin{crllr} \label[crllr]{cor: d and rho commute}
	The permitted operators satisfy
	\begin{align} \label{eq: dpd is d}
		\rhoz^{k + 1} d^k = d^k \rhoz^k & = 0,   &
		d^k p_{k + 1} d^k               & = d^k.
	\end{align}
	As a consequence, $d^k p_{k + 1}$ and $p_{k + 1} d^k$ are idempotent, for all $k$, and $p_\bullet$ is a Poincaré operator in the sense of \Cref{def: poincare operator}.
\end{crllr}
\begin{proof}
	First,  $\rhoz^{k + 1} d^k = 0$ because $\gamz d^k = 0$. Now \Cref{thm: poincare operator} gives us
	\begin{align*}
		0
		= \rhoz^{k + 1} d^k
		= (I - d^k p_{k + 1} - p_{k + 2} d^{k + 1})d^k
		= d^k - d^k p_{k + 1} d^k
		= d^k (I - d^{k - 1} p_k - p_{k + 1} d^k )
		= d^k \rhoz^k.
	\end{align*}
	The second identity is stated after the second equality. The idempotency follows by a left or right application of $p_{k + 1}$.
\end{proof}

\begin{crllr} \label[crllr]{cor: p complex}
	The spaces $V^\bullet$ and operators $p_\bullet$ from \Cref{def: def p} form a Hilbert complex $(\V^\bullet, p_\bullet)$ of negative grade, i.e.
	\begin{align}
		p_{k - 1} p_k = 0
	\end{align}
	for each $k$. Moreover, its cohomology classes have the same dimensions as the cohomology classes of $(\V^\bullet, d^\bullet)$.
\end{crllr}
\begin{proof}
	By \eqref{eq: surjectivity}, we have $\Ran(p_{k + 1}) = \Vb^k$. For $\ub \in \Vb^k$, we have $p_k \ub = 0$ because $\pio \ub = 0$, and thus $\ub \in \Ker(p_k)$. The cohomology class is given by
	\begin{align}
		\Ker(p_k) / \Ran(p_{k + 1})
		= (\Vb^k \oplus \Vz^k) / \Vb^k \simeq \Vz^k
	\end{align}
	and the dimension of $\Vz^k$ equals the dimension of $\Ph^k$ due to \Cref{def: permitting}.
\end{proof}

\begin{rmrk}
	While $p_{k - 1} p_k = 0$ is valid for our construction, we note that any Poincaré operator can be modified to have this property by following \cite[Thm.~5]{vcap2023bounded}.
\end{rmrk}

To obtain the final result of this section, we require an auxiliary lemma. In particular, \Cref{def: permitting} allows us to impose an alternative norm on $\Vb^k$, with which we can formulate a simple, but well-posed, elliptic problem.

\begin{lmm} \label[lmm]{lem: norm on Pb}
	Given a permitting decomposition, then $\| \ub \|_{\Vb^k} \coloneqq \| d^k \ub \|$ is a norm on $\Vb^k$.
\end{lmm}
\begin{proof}
	$\| \cdot \|_{\Vb^k}$ is a semi-norm because $\| \cdot \|$ is a norm and $d^k$ a linear operator. If $\ub \in \Vb^k$ satisfies $d^k \ub = 0$, then $\ub = (\piz \db^k)^{-1} (\piz \db^k) \ub = 0$, showing positive definiteness.
\end{proof}

\begin{crllr} \label[crllr]{lem: wellposed dd}
	For given functional $\accone{f} \in (\Vb^k)^*$, the following problem is well-posed: find $\ub \in \Vb^k$ such that
	\begin{align}
		(d^k \ub, d^k \ub')_W
		             & = \langle \accone{f}, \ub' \rangle, &
		\forall \ub' & \in \Vb^k.
	\end{align}
\end{crllr}
\begin{proof}
	The bilinear form $a(\ub, \ub') \coloneqq (d^k \ub, d^k \ub')_W$ is continuous and coercive in $\| d^k \ub \|$, which is a norm on $\Vb^k$ by \Cref{lem: norm on Pb}. The result now follows by the Lax-Milgram theorem \cite[Thm.~4.1.6]{boffi2013mixed}.
\end{proof}

The results of this section culminate to a new decomposition of the complex, described in the following theorem.
\begin{thrm} \label{thm: decomposition}
	Given a permitting decomposition of $(\V^\bullet, d^\bullet)$, then each $\V^k$ decomposes as
	\begin{align}
		\V^k = \Vb^k \oplus d^{k - 1} \Vb^{k - 1} \oplus \Ph^k.
	\end{align}
\end{thrm}
\begin{proof}
	The inclusion ``$\supseteq$'' is immediate.
	For ``$\subseteq$'', let $u \in \V^k$ and let $\vb_\rho \in \Vb^{k - 1}$ satisfy
	\begin{align} \label{eq: proj of LU}
		(d^{k - 1} \vb_\rho, d^{k - 1} \vb')_W & = (\rhoz^k u, d^{k - 1} \vb')_W, &
		\forall \vb'                           & \in \Vb^{k - 1},
	\end{align}
	which is a well-posed problem by \Cref{lem: wellposed dd}.
	Next, consider $\phi \coloneqq \rhoz^k u - d^{k - 1} \vb_\rho$. \Cref{cor: d and rho commute} gives us $d^k \phi = 0$. Moreover, it implies that for any $v' \in \V^{k - 1}$,
	\begin{align} \label{eq: orthogonal to d}
		(\phi, d^{k - 1} v')_W
		= (\rhoz^k u - d^{k - 1}\vb_\rho, d^{k - 1}p_k d^{k - 1} v')_W
		= 0,
	\end{align}
	because $p_k d^{k - 1} v' \in \Vb^{k - 1}$. We conclude that $\phi \in \Ph^k$.

	Next, let $\ub \coloneqq p_{k + 1} d^k u \in \Vb^k$ and $\vb \coloneqq p_k u + \vb_\rho \in \Vb^{k - 1}$. Using \Cref{thm: poincare operator}, we derive
	\begin{align*}
		\ub + d^{k - 1}\vb + \phi = p_{k + 1} d^k u + d^{k - 1} (p_k u + \vb_\rho) + (\rhoz^k u - d^{k - 1} \vb_\rho) = (d^{k - 1} p_k + p_{k + 1} d^k + \rhoz^k) u = u.
	\end{align*}

	It remains to show that the subspaces are linearly independent. We start by considering $\Vb^k \cap d^{k - 1} \Vb^{k - 1}$. Let $\ub \in \Vb^k$ satisfy $\ub = d^{k - 1} \vb$ for some $\vb \in \Vb^{k - 1}$. Then $\| \ub \|_{\Vb^k} = \| d^k \ub \| = \| d^k d^{k - 1} \vb \| = 0$, so $\ub = 0$ by \Cref{lem: norm on Pb}. The fact that $d^{k - 1} \Vb^{k - 1} \cap \Ph^k = \{0\}$ follows by orthogonality. Finally, if $\ub \in \Vb^k \cap \Ph^k$, then again $d^k \ub = 0$ and thus $\ub = 0$ by \Cref{lem: norm on Pb}.
\end{proof}

\begin{rmrk}
	We emphasize that the decomposition in \Cref{thm: decomposition} is not the same as the permitting decomposition $\V^k = \Vb^k \oplus \Vo^k \oplus \Vz^k$ since $\Vo^k$ is not equal to $d^{k - 1} \Vb^{k - 1}$ in general. Interestingly, the definitions do coincide in the special case of the Hodge decomposition, cf.~\Cref{ex: Hodge}.
\end{rmrk}

\subsection{Poincaré constants}
\label{sub: Poincare constants}

We now briefly turn our focus to the norm from \Cref{lem: norm on Pb}. We will show that this norm is stronger than the $L^2(\Omega)$ norm, with a constant defined as follows.

\begin{dfntn} \label[dfntn]{def: Poincare constant}
	Let the Poincaré constant on $\Vb^k$ be given by
	\begin{align}
		\cb_k \coloneqq \sup_{\ub \in \Vb^k} \frac{\| \ub \|}{\| d^k \ub \|}.
	\end{align}
\end{dfntn}

The constant $\cb_k$ plays an interesting role with respect to the full spaces $\V^k$. In particular, it provides an upper bound for the Poincaré constant on $\V^k$ and, in turn, can be used to bound the inf-sup constant for the pair $\V^{k - 1} \times \V^k$ from below.

\begin{lmm} \label[lmm]{lem: bound on Poincare}
	The Poincaré constant on $\Vb^k$ forms an upper bound on the Poincaré constant on $\V^k$ (cf.~\cite[Sec.~4.2.3]{arnold2018finite}). In particular, it holds that
	\begin{align}
		\| u \| & \le \cb_k \| d^k u \|,
	\end{align}
	for all $u \in \V^k$ with $u \perp \Ker(d^k)$.
\end{lmm}
\begin{proof}
	Let $u \in \V^k$ with $u \perp \Ker(d^k)$. We then use \Cref{thm: poincare operator} and \Cref{cor: d and rho commute} to obtain
	\begin{align*}
		\| u \|^2
		= (u, d^{k - 1} p_k u + p_{k + 1} d^k u + \rhoz^k u)_W
		= (u, p_{k + 1} d^k u)_W
		\le \| u \| \| p_{k + 1} d^k u \|.
	\end{align*}
	Using this bound, we proceed as follows
	\begin{align} \label{eq: calculation ck}
		\sup_{\substack{u \in \V^k                                            \\ u \perp \Ker(d^k)}} \frac{\| u \|}{\| d^k u \|}
		 & \le \sup_{\substack{u \in \V^k                                     \\ u \perp \Ker(d^k)}} \frac{\| p_{k + 1} d^k u \|}{\| d^k u \|} \nonumber \\
		 & \le \sup_{u' \in \V^k} \frac{\| p_{k + 1} d^k u' \|}{\| d^k u' \|}
		= \sup_{u' \in \V^k} \frac{\| p_{k + 1} d^k u' \|}{\| d^k p_{k + 1} d^k u' \|}
		= \sup_{\ub \in \Vb^k} \frac{\| \ub \|}{\| d^k \ub \|}
		= \cb_k.
	\end{align}
	Where we used \eqref{eq: dpd is d} and the fact that $p_{k + 1} d^k$ is identity on $\Vb^k$, and thereby surjective.
\end{proof}
\begin{rmrk}
	From the calculation \eqref{eq: calculation ck}, we moreover see that $\cb_k$ forms a lower bound on the continuity constant of $p_{k + 1}$. Namely,
	\begin{align*}
		\cb_k
		= \sup_{u' \in \V^k} \frac{\| p_{k + 1} d^k u' \|}{\| d^k u' \|}
		\le \sup_{u' \in \V^k} \frac{\| p_{k + 1} d^k u' \|_{V^k}}{\| d^k u' \|_{V^{k + 1}}}
		\le \sup_{t \in \V^{k + 1}} \frac{\| p t \|_{V^k}}{\| t \|_{V^{k + 1}}}
		\eqqcolon \| p_{k + 1} \|,
	\end{align*}
	with the $V^k$-norm defined in \eqref{eq: V norm}.
	This indicates that $\cb_k$ forms a sharper estimate than $\| p_{k + 1} \|$ for the Poincaré constant on $\V^k$.
\end{rmrk}

\begin{crllr}
	The inf-sup constant for the pair $\V^{k - 1} \times \V^k$ is bounded from below as
	\begin{align}
		\inf_{\substack{u \in \Ran(d^{k - 1})}}
		\sup_{v \in \V^{k - 1}}
		\frac{(d^{k - 1} v, u)_W}{\| v \|_{V^{k - 1}} \| u \|_{V^k}}
		\ge
		\frac{1}{\sqrt{\cb_{k - 1}^2 + 1}}.
	\end{align}
\end{crllr}
\begin{proof}
	Let $u = d^{k - 1} v'$ with $v' \in \V^{k - 1}$ and $v' \perp \Ker(d^{k - 1})$. Then, if we set $v = v'$:
	\begin{align*}
		\inf_{\substack{v' \in \V^{k - 1} \\ v' \perp \Ker(d^{k - 1})}}
		\sup_{v \in \V^{k - 1}}
		\frac{(d^{k - 1} v, d^{k - 1} v')_W}{\| v \|_{V^{k - 1}} \| d^{k - 1} v' \|_{V^k}}
		 & \ge
		\inf_{\substack{v' \in \V^{k - 1} \\ v' \perp \Ker(d^{k - 1})}}
		\frac{ \| d^{k - 1} v' \|}{\| v' \|_{V^{k - 1}}}.
	\end{align*}

	\Cref{lem: bound on Poincare} provides the bound $\| v' \|_{V^{k - 1}}^2 \le (\cb_{k - 1}^2 + 1) \| d^{k - 1} v' \|^2$, which concludes the proof.
\end{proof}

\subsection{An explicit basis for the harmonic forms}

In the proof of \Cref{thm: decomposition}, we used a projection to obtain $\phi \in \Ph^k$. This same procedure can be used to construct an explicit basis for the harmonic forms.

\begin{lmm} \label[lmm]{lem: basis harmonics}
	Let $\{\uz_i\}_i$ be a basis of $\Vz^k$. For each $\uz_i$, let $\vb_i \in \Vb^{k - 1}$ satisfy
	\begin{align} \label{eq: proj of uz}
		(d^{k - 1} \vb_i, d^{k - 1} \vb')_W & = (\rhoz^k \uz_i, d^{k - 1} \vb')_W, &
		\forall \vb'                        & \in \Vb^{k - 1}.
	\end{align}
	Then $\{\rhoz^k \uz_i - d^{k - 1} \vb_i\}_i$ is a basis of $\Ph^k$.
\end{lmm}
\begin{proof}
	The fact that $\rhoz^k \uz_i - d^{k - 1} \vb_i \in \Ph^k$ was shown in \Cref{thm: decomposition}, using \Cref{cor: d and rho commute} and \eqref{eq: orthogonal to d}. Moreover, $\dim \Vz^k = \dim \Ph^k$ by \Cref{def: permitting}. It remains to show that they are linearly independent. Consider $\uz \in \Vz^k$ with $\rhoz^k \uz - d^{k - 1} \vb = 0$. Applying $p_k$, \Cref{cor: d and rho commute} and \Cref{cor: p complex} give us
	\begin{align}
		0
		= p_k (\rhoz^k \uz - d^{k - 1} \vb)
		= p_k (I - p_{k + 1} d^k) \uz - p_k d^{k - 1} \vb
		= - \vb
	\end{align}
	from which we conclude that $\rhoz^k \uz = 0$. Using the definition of $\rhoz^k$, this implies that $\uz = p_{k + 1} d^k \uz$. However, $\uz \in \Vz^k$ and $p_{k + 1} d^k \uz \in \Vb^k$, so $\uz = 0$. We conclude that the operation $\uz_i \mapsto \rhoz^k \uz_i - d^{k - 1} \vb_i$ is injective, and thus the range set is linearly independent.
\end{proof}

\section{A direct solver for the Hodge-Laplace problem}
\label{sec: HL solver}

In this section, we focus on the discretization of the Hodge-Laplace problem \eqref{eq: HL problem cont}. For ease of reference, we restate the problem: for given functionals $g \in (\V^{k - 1})^*$ and $f \in (\V^{k})^*$, find $(v, u, \phi) \in \V^{k - 1} \times \V^k \times \Ph^k$ such that
\begin{subequations}  \label{eq: HL problem}
	\begin{align}
		(v, v')_W - (u, d^{k - 1} v')_W                        & = \langle g, v' \rangle,            &
		\forall v'                                             & \in \V^{k - 1}, \label{eq: HL eq 1}   \\
		(d^{k - 1} v, u')_W + (d^k u, d^k u')_W - (\phi, u')_W & = \langle f, u' \rangle,            &
		\forall u'                                             & \in \V^k, \label{eq: HL eq 2}         \\
		(u, \phi')_W                                           & = \langle e, \phi' \rangle,         &
		\forall \phi'                                          & \in \Ph^k. \label{eq: HL eq 3}
	\end{align}
\end{subequations}

We now use the results from the previous section to create a practical method that solves \eqref{eq: HL problem}.
In particular, \Cref{thm: decomposition} allows us to decompose the test and trial functions in \eqref{eq: HL problem} in components from the subspaces $\Vb^\bullet$, $d^\bullet \Vb^\bullet$, and $\Ph^\bullet$. This effectively splits the Hodge-Laplace problem into a sequence of seven symmetric positive definite problems:
\begin{subequations} \label{eqs: seven probs}
	\begin{enumerate}
		\item Find $\vb \in \Vb^{k - 1}$ such that
		      \begin{align}
			      (d^{k - 1} \vb, d^{k - 1} \vb')_W
			                                 & = \langle f, d^{k - 1} \vb' \rangle,                           &
			      \forall \vb'               & \in \Vb^{k - 1}. \label{eq: prob 1}
			      \intertext{
				      \item With $\vb$ given, find $\wb_v \in \Vb^{k - 2}$ such that
			      }
			      (d^{k - 2} \wb_v, d^{k - 2} \wb')_W
			                                 & = \langle g, d^{k - 2} \wb' \rangle - (\vb, d^{k - 2} \wb')_W, &
			      \forall \wb'               & \in \Vb^{k - 2}. \label{eq: prob 2}
			      \intertext{
				      \item Find $\varphi \in \Ph^{k - 1}$ such that
			      }
			      (\varphi, \varphi')_W
			                                 & = \langle g, \varphi' \rangle - (\vb, \varphi')_W,             &
			      \forall \varphi'           & \in \Ph^{k - 1},  \label{eq: prob cohom v}
			      \intertext{
				      \item Find $\phi_f \in \Ph^k$ such that
			      }
			      (\phi_f, \phi')_W
			                                 & = -\langle f, \phi' \rangle,                                   &
			      \forall \phi'              & \in \Ph^k. \label{eq: prob cohom f}
			      \intertext{
				      \item Find $\ub \in \Vb^{k}$ such that
			      }
			      (d^k \ub, d^k \ub')_W
			                                 & = \langle f, \ub' \rangle + ( \phi_f - d^{k - 1} \vb, \ub')_W, &
			      \forall \ub'               & \in \Vb^k. \label{eq: prob 3}
			      \intertext{
				      \item Given the solutions to \eqref{eq: prob 1}--\eqref{eq: prob 3}, find $\vb_u \in \Vb^{k - 1}$ such that
			      }
			      (d^{k - 1} \vb_u, d^{k - 1} \vb')_W
			                                 & = (\vb + d^{k - 2} \wb_v + \varphi, \vb')_W
			      - (\ub, d^{k - 1} \vb')_W
			      - \langle g, \vb' \rangle, &
			      \forall \vb'               & \in \Vb^{k - 1}. \label{eq: prob 4}
			      \intertext{
				      \item Find $\phi_u \in \Ph^k$ such that
			      }
			      (\phi_u, \phi')_W
			                                 & = \langle e, \phi' \rangle - (\ub, \phi')_W,                   &
			      \forall \phi'              & \in \Ph^k, \label{eq: prob cohom u}
		      \end{align}
	\end{enumerate}

	Finally, we collect the solutions from the subproblems and set
	\begin{align} \label{eq: solution}
		v    & \coloneqq \vb + d^{k - 2} \wb_v + \varphi, &
		u    & \coloneqq \ub + d^{k - 1} \vb_u + \phi_u,  &
		\phi & \coloneqq \phi_f
	\end{align}
\end{subequations}

We emphasize that each of the subproblems is well-posed by \Cref{lem: wellposed dd} or Lax-Milgram. In turn, \eqref{eqs: seven probs} is a direct solver for \eqref{eq: HL problem}. It confirms existence of a solution, as shown in the following lemma.

\begin{lmm} \label[lmm]{lem: solution found}
	The triplet $(v, u, \phi) \in \V^{k - 1} \times \V^k \times \Ph^k$ from \eqref{eq: solution} is the solution to the Hodge-Laplace problem \eqref{eq: HL problem}.
\end{lmm}
\begin{proof}
	First, \eqref{eq: prob 2}, \eqref{eq: prob cohom v}, \eqref{eq: prob 4}, and the orthogonality $d^{k - 2} \Vb^{k - 2} \perp \Ph^{k - 1}$ give us
	\begin{align*}
		(v, d^{k - 2} \wb')_W
		                 & = \langle g, d^{k - 2} \wb' \rangle,                &
		\forall \wb'     & \in \Vb^{k - 2},                                      \\
		(v, \varphi')_W
		                 & = \langle g, \varphi' \rangle,                      &
		\forall \varphi' & \in \Ph^{k - 1},                                      \\
		(v, \vb')_W
		                 & = \langle g, \vb' \rangle +  (u, d^{k - 1} \vb')_W, &
		\forall \vb'     & \in \Vb^{k - 1}.
	\end{align*}
	Next, \Cref{thm: decomposition} allows us to decompose the test function $v' \in \V^{k - 1}$ as $v' = \vb' + d^{k - 2} \wb' + \varphi'$ with $(\vb', \wb', \varphi') \in \Vb^{k - 1} \times \Vb^{k - 2} \times \Ph^{k - 1}$. Since $d^{k - 1} d^{k - 2} \wb' = 0$ and $d^{k - 1} \varphi' = 0$, we obtain
	\begin{align*}
		(v, v')_W - (u, d^{k - 1} v')_W
		 & = (v, \vb' + d^{k - 2} \wb' + \varphi')_W - (u, d^{k - 1} \vb')_W                  \\
		 & = \langle g, \vb' + d^{k - 2} \wb' + \varphi' \rangle
		= \langle g, v'\rangle,
		 & \forall v'                                                        & \in \V^{k - 1}
	\end{align*}

	Similarly, \eqref{eq: prob 1}, \eqref{eq: prob cohom f}, and \eqref{eq: prob 3} give us
	\begin{align*}
		(d^{k - 1} v, d^{k - 1} \vb')_W
		              & = \langle f, d^{k - 1} \vb' \rangle, &
		\forall \vb'  & \in \Vb^{k - 1},                       \\
		-(\phi, \phi')_W
		              & = \langle f, \phi' \rangle,          &
		\forall \phi' & \in \Ph^k,                             \\
		(d^{k - 1} v, \ub')_W
		+ (d^k u, d^k \ub')_W
		- (\phi, \ub')_W
		              & = \langle f, \ub' \rangle,           &
		\forall \ub'  & \in \Vb^k.
	\end{align*}
	Thus, if we decompose $u' \in \V^k$ as $u' = \ub' + d^{k - 1} \vb' + \phi'$ with $(\ub', \vb', \phi') \in \Vb^k \times \Vb^{k - 1} \times \Ph^k$, we obtain
	\begin{align*}
		(d^{k - 1} v, u')_W
		+ (d^k u, d^k u')_W
		- (\phi, u')_W
		 & =
		(d^{k - 1} v, \ub' + d^{k - 1} \vb')_W
		+ (d^k u, d^k \ub')_W
		- (\phi, \ub' + \phi')_W                                          \\
		 & = \langle f, \ub' + d^{k - 1} \vb'  + \phi' \rangle
		= \langle f, u' \rangle,
		 & \forall u'                                          & \in \V^k
	\end{align*}

	Finally, \eqref{eq: prob cohom u} ensures that \eqref{eq: HL eq 3} is satisfied.
	In conclusion, $(v, u, \phi)$ is the solution to \eqref{eq: HL problem}.
\end{proof}

Observe that if each $\V^k$ is finite-dimensional, then each of the subproblems in \eqref{eqs: seven probs} is significantly smaller than the original problem \eqref{eq: HL problem}. \Cref{thm: decomposition} defines a different basis for the spaces $\V^\bullet$ and \eqref{eqs: seven probs} is simply \eqref{eq: HL problem} expressed in that basis. This observation is emphasized in the next lemma, for which we introduce the shorthand notation
\begin{align}
	n_k   & \coloneqq \dim (\V^k),  &
	\nb_k & \coloneqq \dim (\Vb^k), &
	\no_k & \coloneqq \dim (\Vo^k), &
	\nz_k & \coloneqq \dim (\Ph^k).
\end{align}

\begin{lmm} \label[lmm]{lem: dims match}
	If $(\V^\bullet, d^\bullet)$ is a finite-dimensional complex, then the dimension of the original Hodge-Laplace problem \eqref{eq: HL problem} is equal to the sum of dimensions of the subproblems \eqref{eqs: seven probs}.
\end{lmm}
\begin{proof}
	The operator $\pio \db: \Vb^k \to \Vo^{k + 1}$ is invertible by \Cref{def: permitting}, so the dimensions of its domain and range coincide, i.e.~$\nb_k = \no_{k + 1}$ for all $k$. In turn, we obtain $n_{k - 1} = \nb_{k - 1} + \no_{k - 1} + \nz_{k - 1} = \nb_{k - 1} + \nb_{k - 2} + \nz_{k - 1}$ and similarly $n_k = \nb_{k} + \nb_{k - 1} + \nz_{k}$.
\end{proof}

\begin{rmrk}
	In practical applications, $v \in \V^{k - 1}$ may be the variable of primary interest. The calculation from \Cref{lem: dims match} shows that we obtain $v$ by solving symmetric positive definite problems of total size $n_{k - 1}$ instead of the original saddle point system of size $n_{k - 1} + n_{k}$.

	For example, Darcy flow in porous media can be modeled by \eqref{eq: HL problem} with $k = n$, in which $v$ denotes the fluid flux, $u$ the pressure, and $d$ the divergence. If the flux is the important variable, e.g. when modeling contaminant transport, then the procedure \eqref{eqs: seven probs} is particularly attractive because it can be obtained by only solving subproblems \eqref{eq: prob 1}, \eqref{eq: prob 2}, and \eqref{eq: prob cohom v}. We consider this case in more detail in the next remark.
\end{rmrk}

\begin{rmrk}[The Poisson problem in mixed form] \label[rmrk]{rem: mixed Darcy}
	Let $\Omega$ be topologically trivial and consider the Hodge-Laplace problem with $k = n$.
	Since $d^n = 0$ and $d^{n - 1}$ is surjective, the procedure from \eqref{eqs: seven probs} simplifies to the following three problems.
	\begin{subequations} \label{eqs: three step kn}
		\begin{itemize}
			\item Find $\vb \in \Pb^{n - 1}$ such that
			      \begin{align}
				      (\divi \vb, u')_\Omega     & = \langle f, u' \rangle,                                    &
				      \forall u'                 & \in \PL^n.  \label{eq: k=n triangle1}
				      \intertext{
					      \item Find $\wb_v \in \Pb^{n - 2}$ such that
				      }
				      (\curl \wb_v, \curl \wb')_\Omega
				                                 & = \langle g, \curl \wb' \rangle - (\vb, \curl \wb')_\Omega, &
				      \forall \wb'               & \in \Pb^{n - 2}. \label{eq: k = n eq2}
				      \intertext{
					      \item Find $u \in \PL^n$ such that
				      }
				      (u, \divi \vb')_\Omega
				                                 & = (\vb + \curl \wb_v, \vb')_\Omega
				      - \langle g, \vb' \rangle, &
				      \forall \vb'               & \in \Pb^{n - 1}. \label{eq: k=n triangle2}
			      \end{align}
		\end{itemize}
	\end{subequations}
	Finally, we set $v = \vb + d \wb_v$.

	In the case of \Cref{ex: spanning tree decomposition}, \eqref{eq: k=n triangle1} and \eqref{eq: k=n triangle2} are easy to solve due to the underlying tree structure, as observed in \cite{boon2025deep}. These can therefore be solved in linear time with respect to the number of elements. The remaining subproblem \eqref{eq: k = n eq2} is a symmetric, positive definite system of size $\nb_{n - 2}$, which corresponds to the number of faces \emph{minus} the number of cells.

	Note that \eqref{eq: k = n eq2} is exactly \eqref{eq: HL eq 1} posed in a divergence-free subspace of $\PL^{n - 1}$. This relates directly to the divergence-free finite element spaces proposed in \cite{alonso2017graphs,devloo2022efficient}, which are examples of the form $d^{n - 2} \Pb^{n - 2}$. In particular, they correspond to the curl of a subset of the $H(\curl)$-conforming elements $\PL^{n - 2}$, filtered by a spanning tree.
\end{rmrk}

As shown in \Cref{lem: solution found}, the true solution is found by \eqref{eqs: seven probs} if exact solvers are used. However, it may be advantageous to use inexact, iterative solvers in some of the steps. As shown in \Cref{lem: wellposed dd}, the problems \eqref{eqs: seven probs} are symmetric and positive definite, and thus amenable to efficient solvers such as the Conjugate Gradient method. We highlight an important implication of such solvers in the following lemma.

\begin{lmm}
	If inexact solvers are used in \eqref{eqs: seven probs}, then the error in \eqref{eq: HL eq 2} depends only on the accuracy at which problems \eqref{eq: prob 1}, \eqref{eq: prob cohom f}, and \eqref{eq: prob 3} are solved.
\end{lmm}
\begin{proof}
	Let $(v_*, u_*) \in \V^{k - 1} \times \V^k$ be the solution obtained using inexact solvers with $v_* = \vb_* + d^{k - 2} \wb_{v,*} + \varphi_*$ and $u_* = \ub_* + d^{k - 1} \vb_{u,*} + \phi_*$. The error equation for \eqref{eq: HL eq 2} becomes
	\begin{align*}
		(d^{k - 1} v_*, u')_\Omega + (du_*, du')_\Omega - \langle f, u' \rangle
		 & = (d^{k - 1} \vb_*, u')_\Omega + (d^k \ub_*, du')_\Omega - \langle f, u' \rangle.
	\end{align*}
	In other words, the error only depends on the components $(\vb_*, \ub_*)$, which are obtained from \eqref{eq: prob 1}, \eqref{eq: prob cohom f}, and \eqref{eq: prob 3}.
\end{proof}

\section{Auxiliary space preconditioning}
\label{sec: aux preconditioning}

The solver proposed in \Cref{sec: HL solver} depends on the particular structure of the Hodge-Laplace problem and uses it to its advantage. In this section, we aim to use the results from \Cref{sec: the Poincare operator} to form a solver that is applicable to a more general class of problems.

For that, we use the framework of norm-equivalent preconditioning \cite{mardal2011preconditioning}, which employs the Riesz representation operator as the canonical preconditioner for (mixed) finite element problems.
In short, if the problem is well-posed in a certain, parameter-weighted norm, then that norm determines a suitable preconditioner.
Such a norm is often a weighted variant of the $H \Lambda^k$ norm from \eqref{eq: H norm} so let us consider
\begin{align} \label{eq: weighted U norm}
	\| u \|_U^2 \coloneqq \| \alpha u \|_\Omega^2 + \| d^k u \|_\Omega^2.
\end{align}
in which $0 < \alpha \le 1$ is a fixed constant.
Following \cite{mardal2011preconditioning}, if our problem is well-posed in this norm, then the canonical preconditioner involves solving the following projection problem: given $f \in (\V^k)^*$, find $u \in P\Lambda^k$ such that
\begin{align} \label{eq: projection problem}
	(\alpha^2 u, u')_\Omega
	+ (d^k u, d^k u')_\Omega
	           & = \langle f, u' \rangle, &
	\forall u' & \in \V^k.
\end{align}

However, problem \eqref{eq: projection problem} may be challenging to solve numerically if $\alpha$ is small, because the term $(du, du')_\Omega$ has a large kernel for $0 < k < n$.
To handle this, we follow the framework of \cite{hiptmair2007nodal,xu1996auxiliary} and construct an auxiliary space preconditioner based on the decomposition from \Cref{thm: decomposition}. In particular, we define the auxiliary spaces $U_j$ and transfer operators $\Pi_j: U_j \to \V^k$ as
\begin{subequations} \label{eqs: transfers}
	\begin{align}
		U_1              & \coloneqq \Vb^k,                             &
		\| \ub \|_{U_1}  & \coloneqq \| d^k \ub \|_\Omega,              &
		\Pi_1 \ub        & \coloneqq \ub,                                 \\ 
		U_2              & \coloneqq \Vb^{k - 1},                       &
		\| \vb \|_{U_2}  & \coloneqq \| \alpha d^{k - 1} \vb \|_\Omega, &
		\Pi_2 \vb        & \coloneqq d^{k - 1} \vb,                       \\
		U_3              & \coloneqq \Ph^k,                             &
		\| \phi \|_{U_3} & \coloneqq \| \alpha \phi \|_\Omega,          &
		\Pi_3 \phi       & \coloneqq \phi.                                
	\end{align}
\end{subequations}

Next, we show two stability results, concerning first the transfer operators and second, the decomposition.

\begin{lmm} \label[lmm]{lem: stable transfer}
	The transfer operators $\Pi_j$ from \eqref{eqs: transfers} are stable. I.e.~for each $j \in \{1, 2\}$, a bounded $c_j \in \mathbb{R}$ exists such that
	\begin{align}
		\| \Pi_j u_j \|_U & \le c_j \| u_j \|_{U_j}, &
		\forall u_j       & \in U_j,
	\end{align}
	where $c_j$ may depend on $\cb_k$ from \Cref{def: Poincare constant}, but not on $\alpha$.
\end{lmm}
\begin{proof}
	The result follows immediately from \Cref{lem: bound on Poincare}, the upper bound on $\alpha$, and the identities $d^{k + 1} d^k = 0$ and $d^k \phi = 0$:
	\begin{subequations}
		\begin{align}
			\| \alpha \ub \|_\Omega^2 + \| d^k \ub \|_\Omega^2
			                                          & \le ((\alpha \cb_k)^2 + 1) \| d^k \ub \|_\Omega^2
			\le (\cb_k^2 + 1) \| d^k \ub \|_\Omega^2, &
			\forall \ub                               & \in \Vb^k, \label{eq: transfer1 stable}             \\
			\| \alpha d^{k - 1} \vb \|_\Omega^2 + \| d^k (d^{k - 1} \vb) \|_\Omega^2
			                                          & = \| \alpha d^{k - 1} \vb \|_\Omega^2,            &
			\forall \vb                               & \in \Vb^{k - 1},                                    \\
			\| \alpha \phi \|_\Omega^2 + \| d^k \phi \|_\Omega^2
			                                          & = \| \alpha \phi \|_\Omega^2,                     &
			\forall \phi                              & \in \Ph^k.
		\end{align}
	\end{subequations}
\end{proof}

\begin{lmm} \label[lmm]{lem: stable decomposition}
	The decomposition from \Cref{thm: decomposition} is stable in the norms from \eqref{eq: weighted U norm} and \eqref{eqs: transfers}. In particular, if $u \in \V^k$ decomposes as $u = \ub + d^{k - 1} \vb + \phi$ with $(\ub, \vb, \phi) \in \Vb^k \times \Vb^{k -1} \times \Ph^k$, then
	\begin{align}
		\| d^k \ub \|_\Omega + \| \alpha d^{k - 1} \vb \|_\Omega
		+ \| \alpha \phi \|_\Omega
		\le (1 + 4 \| \rhoz^k \| )\| \alpha u \|_\Omega + (1 + \cb_k) \| d^k u \|_\Omega.
	\end{align}
\end{lmm}
\begin{proof}
	Recall from \Cref{thm: decomposition} that the decomposition is uniquely given by
	$\ub \coloneqq p_{k + 1} d^k u \in \Vb^k$, $\vb \coloneqq p_k u + \vb_\rho \in \Vb^{k - 1}$, and $\phi \coloneqq \rhoz^k u - d^{k - 1} \vb_\rho$, with $\vb_\rho$ given by \eqref{eq: proj of LU}.
	Using the identity \eqref{eq: dpd is d}, we derive
	\begin{align} \label{eq: stability 1}
		\| d^k \ub \|_\Omega = \| d^k p_{k + 1} d^k u \|_\Omega = \| d^k u \|_\Omega.
	\end{align}

	Next, \eqref{eq: proj of LU} provides the bound $\| d^{k - 1} \vb_\rho \|_\Omega \le \| \rhoz^k u \|_\Omega$. We use this to bound $\phi$ as follows
	\begin{align} \label{eq: stability 1.5}
		\| \alpha \phi \|_\Omega
		= \| \alpha (\rhoz^k u - d^{k - 1} \vb_\rho) \|_\Omega
		\le 2 \| \alpha \rhoz^k u \|_\Omega
		\le 2 \| \rhoz^k \| \| \alpha u \|_\Omega.
	\end{align}

	Finally, we bound $\vb$ by using the decomposition $d^{k - 1} \vb = u - \ub - \phi$, \Cref{lem: bound on Poincare}, and the upper bound on $\alpha$
	\begin{align} \label{eq: stability 2}
		\| \alpha d^{k - 1} \vb \|_\Omega
		\le \| \alpha u \|_\Omega + \| \alpha \ub \|_\Omega + \| \alpha \phi \|_\Omega
		 & \le (1 + 2 \| \rhoz^k \|) \| \alpha u \|_\Omega + \alpha \cb_k \| d^k \ub \|_\Omega \nonumber \\
		 & \le (1 + 2 \| \rhoz^k \|) \| \alpha u \|_\Omega + \cb_k \| d^k u \|_\Omega.
	\end{align}
	Together, \eqref{eq: stability 1}, \eqref{eq: stability 1.5}, and \eqref{eq: stability 2} yield the result.
\end{proof}

Let $\mathsf{\accone{\Pi}_k} \in \mathbb{R}^{\accone{n}_k \times n_k}$ be the matrix representation of the restriction operator $\pib$ on $\V^k$. Similarly, let $\mathsf{\accone{D}_{k - 1}} \in \mathbb{R}^{n_k \times \nb_{k - 1}}$ represent the differential $\db$ on $\Vb^{k - 1}$. Let $\mathsf{\accone{A}_k} \in \mathbb{R}^{\accone{n}_k \times \accone{n}_k}$ be the Gram matrix for the norm from \Cref{lem: norm on Pb}. In other words,
\begin{align*}
	\| d^k \ub \|_\Omega^2 & = \mathsf{\ub^T \accone{A}_k \ub}, &
	\forall \ub            & \in \Vb^k,
\end{align*}
in which $\mathsf{\ub} \in \mathbb{R}^{n_k}$ is the vector representation of $\ub$. Similarly, let $\mathsf{\accthr{\Pi}_k} \in \mathbb{R}^{\nz_k \times n_k}$ be the matrix that has the basis functions of $\Ph^k$ as its columns, and let $\accthr{\mathsf{M}}_k \in \mathbb{R}^{\nz_k \times \nz_k}$ denote the mass matrix of $\Ph^k$. Using these ingredients, we now construct our preconditioner.

\begin{dfntn} \label[dfntn]{def: aux space precond}
	The auxiliary space preconditioner for \eqref{eq: projection problem}, based on the decomposition from \Cref{thm: decomposition}, is given by
	\begin{align}
		\mathsf{P \coloneqq \accone{\Pi}_k^T \accone{A}_k^{-1} \accone{\Pi}_k}
		+ \mathsf{\accone{D}_{k - 1}}
		(\alpha^2 \mathsf{\accone{A}_{k - 1})^{-1} \accone{D}_{k - 1}^T}
		+ \mathsf{\accthr{\Pi}_k}
		(\alpha^2 \mathsf{\accthr{M}_{k - 1})^{-1} \accthr{\Pi}_k^T}.
	\end{align}
\end{dfntn}

\begin{lmm} \label[lmm]{lem: robust precond}
	The auxiliary space preconditioner $\mathsf{P}$ from \Cref{def: aux space precond} is robust with respect to $0 < \alpha \le 1$. Moreover, if $\cb_k$ and $\| \rhoz^k \|$ are independent of the mesh size $h$, then $\mathsf{P}$ is robust with respect to $h$.
\end{lmm}
\begin{proof}
	We verify the assumptions of \cite[Thm. 2.2]{hiptmair2007nodal}. First, \Cref{thm: decomposition} provides the decomposition and therefore the surjectivity of the combined transfer operators. The continuity of these operators is shown in \Cref{lem: stable transfer} and \Cref{lem: stable decomposition} provides the stability of the decomposition.

	Next, we note that the bounds in these results are independent of $\alpha$ but depend on $\cb_k$ and $\| \rhoz^k \|$. In turn, \cite[Thm. 2.2]{hiptmair2007nodal} shows that the spectral condition number of the preconditioned system is independent of $\alpha$. Moreover, if $\cb_k$ and $\| \rhoz^k \|$ are independent of $h$, then the preconditioner is robust with respect to $h$ as well.
\end{proof}

\section{Numerical results}
\label{sec: Numerical results}

We now present three numerical experiments that highlight the practical impact of the results derived in the previous sections. \Cref{sub: num HL} concerns the seven-step solution procedure from \Cref{sec: HL solver}, \Cref{sub: num harmonics} presents a basis for the harmonic forms based on \Cref{lem: basis harmonics}, \Cref{sub: num poincare} investigates the Poincaré constant from \Cref{sub: Poincare constants}, and \Cref{sub: num aux} considers the auxiliary space preconditioners from \Cref{sec: aux preconditioning}.

Throughout this section, we define $\V^\bullet \coloneqq \PL^\bullet$ as the Whitney forms from \Cref{def: whitney} and consider the permitting decomposition described in \Cref{ex: spanning tree decomposition}. We construct the spanning trees as follows. Let $\mathcal{T}_{n - 1}^*$ be the spanning tree of the dual graph $\mathcal{G}_{n - 1}^*$, obtained by a breadth-first search algorithm rooted at $m^\infty$, the vertex that corresponds to the ``outside cell''. We note that this construction is similar to \cite[Sec.~7.3]{jiranek2010posteriori}. Moreover, \cite[Sec.~6]{alonso2018graph} showed that a breadth-first is preferred over a depth-first search algorithm. If $n = 2$, we define $\mathcal{E}_1$ as the complement to $\mathcal{F}_1$. An illustration in 2D is given in \Cref{fig: trees}. Finally, for $n = 3$, we define $\mathcal{T}_1$ as the spanning tree of the nodes, using again a breadth-first search, rooted at the node that is closest to the center of the domain.

The code is written in Python 3, uses the open-source package PyGeoN\cite{pygeon}, and is freely available at \url{github.com/wmboon/poincare}.

\subsection{Solving the Hodge-Laplace problem}
\label{sub: num HL}

In this section, we investigate the efficiency of the solution procedure \eqref{eqs: seven probs} proposed in \Cref{sec: HL solver}.
The experimental set-up is as follows. Let $\Omega_h$ be an unstructured simplicial grid of the unit cube consisting of approximately \num{17200} elements. The edges have a mean length of $h \approx 0.08$. We set up the Hodge-Laplace problem with randomly chosen right-hand sides. We then compare the solution times needed for a direct solver for each of the subproblems of \eqref{eqs: seven probs} with the time required to solve the saddle-point system \eqref{eq: HL problem}. For each problem, we report the number of degrees of freedom and the time required by the direct solver. The results are presented in \Cref{tab: solve times}.

The experiment is then repeated for a domain consisting of the unit cube minus a torus. The domain has the Betti numbers $(1, 1, 1, 0)$ and is illustrated in \Cref{fig:harmonic1}. The mesh consists of approximately \num{21400} elements, and the bases for the non-trivial harmonic forms is computed using \Cref{lem: basis harmonics}. We report the solver times in \Cref{tab: solve times nontrivial} and make four observations that relate to both cases.

\begin{table}[ht]
	\caption{Dimensions and solution times for direct solvers of the subproblems in \eqref{eqs: seven probs} compared to the original system \eqref{eq: HL problem}, on the unit cube. Since $\Ph^k = \{0\}$ for all $k > 0$, we omit those rows in the table.
	}
	\label{tab: solve times}
	\centering

	\begin{tabular}{ll|rr|rr|rr}
		\hline
		                &
		                & \multicolumn{2}{c|}{$k = 1$}
		                & \multicolumn{2}{c|}{$k = 2$}
		                & \multicolumn{2}{c}{$k = 3$}
		\\
		Space           & Problem
		                & DoF                          & Time
		                & DoF                          & Time
		                & DoF                          & Time
		\\
		\hline
		                &                              &              &       &              &        &                      \\[-2ex]
		$ \Pb^{k-1}$    & \eqref{eq: prob 1}           & \num{  3741} & 0.05s & \num{ 19035} & 2.32s  & \num{ 17206} & 0.01s \\
		$ \Pb^{k-2}$    & \eqref{eq: prob 2}           & \num{     0} & -     & \num{  3741} & 0.05s  & \num{ 19035} & 2.85s \\
		$ \Ph^{k-1}$    & \eqref{eq: prob cohom v}     & \num{     1} & 0.00s & \num{     0} & -      & \num{     0} & -     \\
		$ \Pb^k$        & \eqref{eq: prob 3}           & \num{ 19035} & 2.36s & \num{ 17206} & 0.01s  & \num{     0} & -     \\
		$ \Pb^{k-1}$    & \eqref{eq: prob 4}           & \num{  3741} & 0.05s & \num{ 19035} & 2.35s  & \num{ 17206} & 0.01s \\
		\hline
		Total           & \eqref{eqs: seven probs}     & \num{ 26518} & 2.46s & \num{ 59017} & 4.73s  & \num{ 53447} & 2.87s \\
		\hline \hline
		Original system & \eqref{eq: HL problem}       & \num{ 26518} & 8.71s & \num{ 59017} & 43.88s & \num{ 53447} & 5.77s
	\end{tabular}
\end{table}

\begin{table}[ht]
	\caption{Dimensions and solution times for direct solvers of the seven subproblems in \eqref{eqs: seven probs} compared to the original system \eqref{eq: HL problem}, on the unit cube with a removed torus.
	}
	\label{tab: solve times nontrivial}
	\centering

	\begin{tabular}{ll|rr|rr|rr}
		\hline
		                &
		                & \multicolumn{2}{c|}{$k = 1$}
		                & \multicolumn{2}{c|}{$k = 2$}
		                & \multicolumn{2}{c}{$k = 3$}
		\\
		Space           & Problem
		                & DoF                          & Time
		                & DoF                          & Time
		                & DoF                          & Time
		\\
		\hline
		                &                              &              &       &              &        &                      \\[-2ex]
		$ \Pb^{k-1}$    & \eqref{eq: prob 1}           & \num{  4707} & 0.09s & \num{ 23921} & 3.30s  & \num{ 21444} & 0.01s \\
		$ \Pb^{k-2}$    & \eqref{eq: prob 2}           & \num{     0} & -     & \num{  4707} & 0.10s  & \num{ 23921} & 3.59s \\
		$ \Ph^{k-1}$    & \eqref{eq: prob cohom v}     & \num{     1} & 0.01s & \num{     1} & 0.00s  & \num{     1} & 0.00s \\
		$ \Ph^k$        & \eqref{eq: prob cohom f}     & \num{     1} & 0.00s & \num{     1} & 0.00s  & \num{     0} & -     \\
		$ \Pb^k$        & \eqref{eq: prob 3}           & \num{ 23921} & 2.84s & \num{ 21444} & 0.01s  & \num{     0} & -     \\
		$ \Pb^{k-1}$    & \eqref{eq: prob 4}           & \num{  4707} & 0.09s & \num{ 23921} & 3.70s  & \num{ 21444} & 0.02s \\
		$ \Ph^k$        & \eqref{eq: prob cohom u}     & \num{     1} & 0.00s & \num{     1} & 0.00s  & \num{     0} & -     \\
		\hline
		Total           & \eqref{eqs: seven probs}     & \num{ 33338} & 3.03s & \num{ 73996} & 7.11s  & \num{ 66810} & 3.62s \\
		\hline \hline
		Original system & \eqref{eq: HL problem}       & \num{ 33338} & 9.61s & \num{ 73996} & 87.51s & \num{ 66810} & 8.31s
	\end{tabular}
\end{table}

First, while it is unsurprising that the total numbers of degrees of freedom coincide (cf.~\Cref{lem: dims match}), it is noteworthy that the problem for $k = 3$ is decomposed into three subproblems of roughly equal size. As noted in \Cref{rem: mixed Darcy}, this case corresponds to the mixed formulation of the Poisson equation that is common in Darcy flow models.

Second, we observe a speed-up factor of at least two for all instances of $k$, without loss of accuracy, which indicates the potential of the approach. The best relative performance is seen in the case of $k = 2$, i.e.~the Hodge-Laplace problem posed on the Raviart-Thomas space, where the speed-up factor exceeds 10 for both set-ups. Interestingly, this is the only problem that requires solving all seven subproblems. The problems posed on $\Vb^{2}$ have an underlying tree structure which allows them to be solved in sequence, from the leaves to the root. We have not explicitly informed the solver of this and implemented exactly \eqref{eqs: seven probs}. Nonetheless, these problems are solved within two hundredths of a second. On the other hand, the problems on $\Vb^{0}$ are simply much smaller than the others because there are far fewer nodes than edges or faces in the grid. These problems are therefore also solved relatively quickly.

Third, we note that each of the (sub)problems can be solved with more efficient solution methods. However, we have opted for a standard, direct linear solver in all cases to keep the comparison relatively fair. The performance of different iterative solvers forms a study in its own right, which we reserve as a topic for future investigation.

Finally, assembly times were not included in this comparison. The assembly of \eqref{eqs: seven probs} is similar to the assembly of \eqref{eq: HL problem} because each term consists of the same building blocks. In particular, for given $k$, we require the mass matrices of $\PL^{k - 1}$, $\PL^k$, $\PL^{k + 1}$, $\Ph^k$, and $\Ph^{k - 1}$. The matrix representations of the differentials, $\mathsf{D}_k$, can directly be obtained from the incidence matrices of the mesh. All that remains is the restriction onto the subspaces $\Vb^k$, which can be done using the same machinery as imposing essential boundary conditions. We remark that one can even apply the solution procedure \eqref{eqs: seven probs} without assembly of the original system \eqref{eq: HL problem}. If direct solvers are used, as they are here, then \Cref{lem: solution found} shows that the same solution is obtained.

\subsection{A basis for the harmonic forms}
\label{sub: num harmonics}

\begin{figure}[ht]
	\centering
	\includegraphics[width=0.9\linewidth]{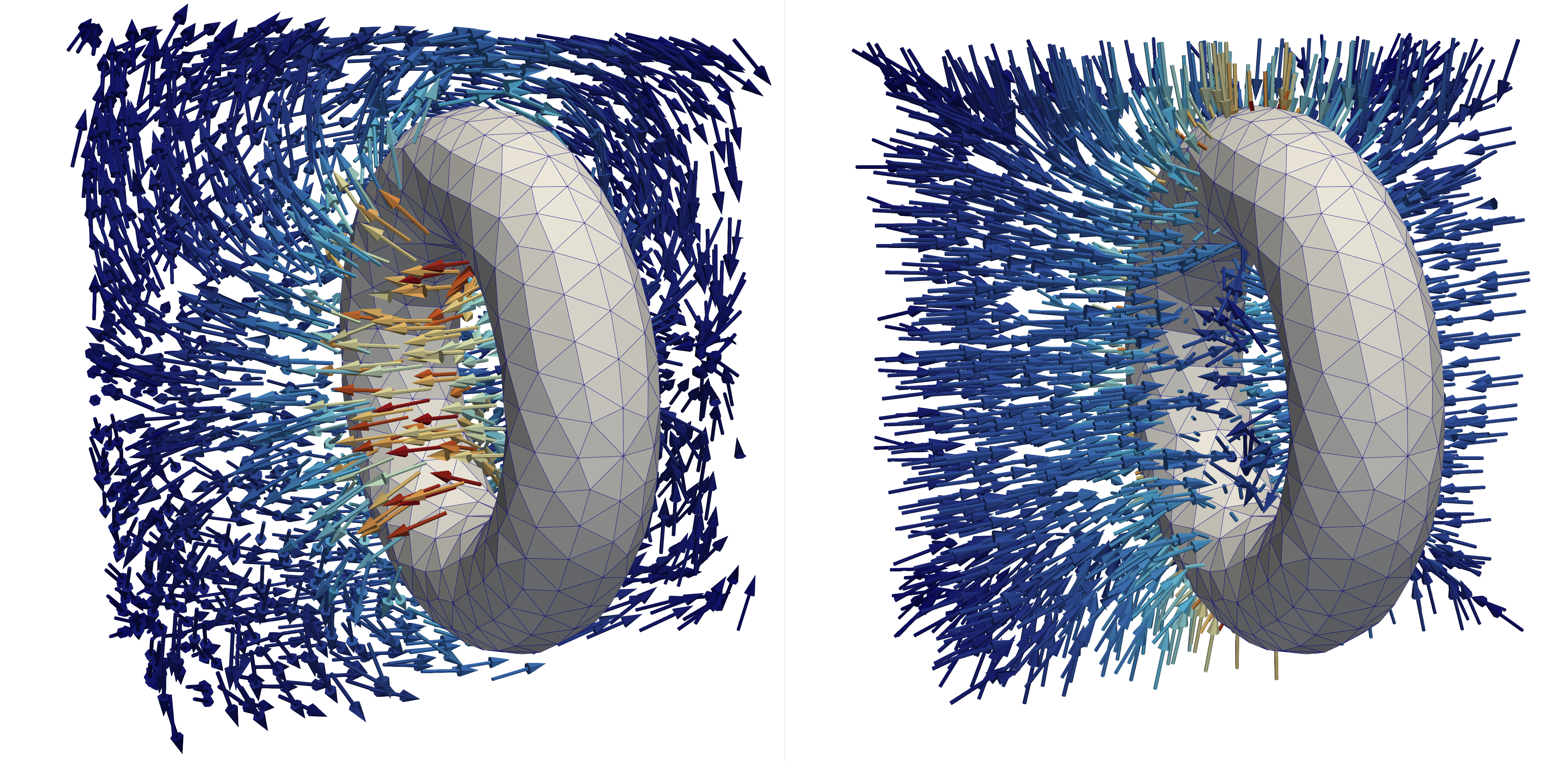}
	\caption{The basis vectors for the harmonic forms in the case of a torus removed from a cube. Only half of the vector fields are shown for visualization purposes. }
	\label{fig:harmonic1}
\end{figure}

The second test case in the previous experiment has non-trivial harmonic 1-forms and 2-forms. Using the technique from \Cref{lem: basis harmonics}, we compute explicit bases for these spaces, which are illustrated in \Cref{fig:harmonic1}.

In this case, the continuous harmonic forms are vector fields that are solenoidal and irrotational. In the discrete setting, the harmonic 1-forms are irrotational and orthogonal to gradients from the discrete nodal space. On the other hand, the discrete harmonic 2-forms are solenoidal and orthogonal to curls. The computed basis functions illustrated in \Cref{fig:harmonic1} confirm the expected behavior mentioned in the first two points of \Cref{lem: permitting}. In particular, the harmonic 1-forms circulate around the torus, whereas the harmonic 2-forms are vector fields that point from the outer to the inner boundary. This interpretation is opposite for the complex with boundary conditions.

The cost of finding these harmonic forms is as follows. Constructing the spanning tree decomposition of \Cref{ex: spanning tree decomposition} is computationally inexpensive. Then $\rhoz^k \uz_i$ contains one application of $p_{k + 1}$ which involves solving a system of size $\no_{k + 1} = \nb_k$. Finally, the projection problem \eqref{eq: proj of uz} involves a symmetric positive definite system of size $\nb_{k - 1} = \no_k$. In total, we therefore only compute for the relevant $\nb_k + \no_k = n_k - \nz_k$ degrees of freedom.

\subsection{Estimating the Poincaré constant}
\label{sub: num poincare}

To compute the constant from \Cref{def: Poincare constant}, we consider the following generalized eigenvalue problem:
find the largest $\lambda \in \mathbb{R}$ such that a $\ub \in \Vb^k$ exists with
\begin{align}
	(\ub, \ub')_\Omega & = \lambda (d^k \ub, d^k \ub')_\Omega, &
	\forall \ub' \in \Vb^k.
\end{align}
It then follows that $\cb_k = \sqrt{\lambda}$. We note that this problem can be handled numerically because the matrices on the left and right-hand sides are both invertible (cf.~\Cref{lem: norm on Pb}). This allows us to compute $\cb_k$ for a series of meshes on the unit square and cube. We present the computed constants in \Cref{tab:constant_estimate}. The mesh used in \Cref{sub: num HL} corresponds to the finest level on the unit cube.

\begin{table}[ht]
	\caption{Computed constants $\cb_k$ (cf.~\Cref{def: Poincare constant}) on the unit square and cube. \Cref{lem: bound on Poincare} shows that $\cb_k$ forms an upper bound on the Poincaré constant of $\PL^k$.}
	\label{tab:constant_estimate}

	\centering

	\begin{tabular}{r|rr||r|rrr}
		\hline
		\multicolumn{3}{c||}{Unit square} &
		\multicolumn{4}{c}{Unit cube}                                                               \\
		\hline
		$h$                               & $k = 0$ & $k = 1$ &
		$h$                               & $k = 0$ & $k = 1$ & $k = 2$                             \\
		\hline
		1.29e-01                          & 0.317   & 0.377   & 5.39e-01 & 0.327  & 1.408  & 0.649  \\
		6.52e-02                          & 0.318   & 0.485   & 3.94e-01 & 0.310  & 0.769  & 0.511  \\
		3.35e-02                          & 0.318   & 0.604   & 3.00e-01 & 0.313  & 1.006  & 0.607  \\
		1.66e-02                          & 0.318   & 0.841   & 2.14e-01 & 0.315  & 1.126  & 0.632  \\
		8.25e-03                          & 0.318   & 1.034   & 1.48e-01 & 0.317  & 1.544  & 0.763  \\
		\hline
		\hline
		Rate                              & 0.000   & -0.368  &          & -0.019 & -0.646 & -0.343
	\end{tabular}
\end{table}

The case $k = 0$ provides the most stable constants because $\Vb^0$ is the nodal space constrained to be zero in one node (cf.~\Cref{ex: spanning tree decomposition}). In turn, $\cb_0$ is exactly the Poincaré constant for the gradient on the nodal space $\PL^0$. For both cases of $n$, the constant approximates $1 / \pi$.

The other constants exhibit a mild increase with respect to $h$. Assuming $\cb_k \approx \mathcal{O}(h^\alpha)$, we compute the rate $\alpha$ based on the three finest meshes and report these in the final row. We observe that the inverse dependency is sub-linear in all cases. We remark that these are the constants and rates that result from our specific choice of spanning trees. It may be possible to obtain better estimates with a different construction, using e.g.~the multi-level trees of \cite{hiptmair2000multilevel}.

\subsection{Preconditioning the projection problem}
\label{sub: num aux}

Let us consider the projection problem \eqref{eq: projection problem} with the preconditioner from \Cref{def: aux space precond}. We note that the only relevant cases are $1 \le k \le n - 1$, leaving three non-trivial cases. Using grids on the unit cube and square, we investigate its performance for a range of $10^{-4} \le \alpha \le 1$ by applying Conjugate Gradient method, preconditioned with the auxiliary space preconditioner, denoted as ASP. We report the number of iterations required to reach a relative residual of $10^{-5}$.

Direct solvers are not feasible for larger systems, so we also consider a variant in which the linear solves in \Cref{def: aux space precond} are replaced by a V-cycle of smoothed aggregation algebraic multigrid. We denote this variant by ASP-AMG. We do apply a direct solver for the matrix $\accone{\mathsf{A}}_{n - 1}$ because these systems are fast to solve, as shown in the first experiment.

The third preconditioner we consider is the smoothed aggregation algebraic multigrid on the original system (AMG). All AMG operators are implemented using the default settings of PyAMG \cite{bell2022pyamg}. This experiment therefore does not capture the full capabilities of AMG with optimized smoothers. Rather, we aim to show that the application of the decomposition can improve the performance of ``out-of-the-box" solvers.

\begin{table}[ht]
	\caption{Performance of the auxiliary space preconditioner of \Cref{def: aux space precond} for the projection problem \eqref{eq: projection problem}. A dash indicates that the conjugate gradient method did not converge within 750 iterations.}
	\label{tab:aux_space}

	\small
	\centering

	\begin{tabular}{r|rrrrr|rrrrr|rrrrr}
		\hline
		                             &
		\multicolumn{5}{c|}{ASP}     &
		\multicolumn{5}{c|}{ASP-AMG} &
		\multicolumn{5}{c}{AMG}                                                                                           \\
		$\log_{10}(\alpha)$          & -4 & -3 & -2 & -1 & 0
		                             & -4 & -3 & -2 & -1 & 0
		                             & -4 & -3 & -2 & -1 & 0                                                              \\
		\hline
		$h$                          &
		\multicolumn{15}{c}{Unit square, $k = 1$}                                                                         \\
		\hline
		1.29e-01                     & 2  & 2  & 3  & 4  & 9  & 9   & 9   & 9   & 9   & 10  & 641 & 407 & 259 & 154 & 72  \\
		6.52e-02                     & 3  & 2  & 4  & 4  & 10 & 10  & 10  & 10  & 10  & 12  & -   & -   & -   & 431 & 165 \\
		3.35e-02                     & 3  & 2  & 4  & 5  & 12 & 11  & 11  & 11  & 11  & 14  & -   & -   & -   & -   & 366 \\
		1.66e-02                     & 3  & 3  & 4  & 6  & 16 & 13  & 13  & 13  & 13  & 18  & -   & -   & -   & -   & -   \\
		8.25e-03                     & 4  & 3  & 4  & 6  & 19 & 16  & 16  & 16  & 16  & 24  & -   & -   & -   & -   & -   \\
		\hline
		$h$                          &
		\multicolumn{15}{c}{Unit cube, $k = 1$}                                                                           \\
		\hline
		9.02e-01                     & 2  & 2  & 3  & 4  & 10 & 9   & 9   & 9   & 9   & 12  & 89  & 54  & 37  & 28  & 15  \\
		5.84e-01                     & 2  & 2  & 4  & 6  & 14 & 18  & 17  & 17  & 17  & 20  & 202 & 144 & 96  & 60  & 32  \\
		3.79e-01                     & 2  & 2  & 4  & 5  & 13 & 25  & 25  & 25  & 24  & 27  & 349 & 268 & 182 & 113 & 51  \\
		2.01e-01                     & 3  & 2  & 4  & 6  & 19 & 66  & 66  & 65  & 65  & 68  & -   & 552 & 362 & 219 & 95  \\
		1.02e-01                     & 3  & 3  & 4  & 7  & 30 & 205 & 204 & 204 & 203 & 209 & -   & -   & -   & 544 & 197 \\
		\hline
		$h$                          &
		\multicolumn{15}{c}{Unit cube, $k = 2$}                                                                           \\
		\hline
		9.02e-01                     & 2  & 2  & 3  & 4  & 9  & 10  & 10  & 9   & 9   & 12  & 170 & 115 & 79  & 54  & 26  \\
		5.84e-01                     & 2  & 2  & 4  & 4  & 12 & 21  & 21  & 19  & 19  & 23  & 392 & 300 & 202 & 134 & 53  \\
		3.79e-01                     & 2  & 2  & 4  & 4  & 10 & 33  & 33  & 31  & 29  & 33  & -   & 553 & 404 & 250 & 90  \\
		2.01e-01                     & 3  & 2  & 4  & 5  & 11 & 99  & 94  & 86  & 79  & 88  & -   & -   & 723 & 462 & 159 \\
		1.02e-01                     & 3  & 2  & 4  & 5  & 15 & 307 & 287 & 261 & 243 & 293 & -   & -   & -   & -   & 334 \\
	\end{tabular}
\end{table}

As shown in \Cref{tab:aux_space}, the ASP with direct solvers is particularly effective for small values of $\alpha$, often letting the conjugate gradient method converge in two iterations. That case is typically challenging to solve since the second term in \eqref{eq: projection problem} dominates, which has a large kernel. Because the preconditioner is based on the decomposition from \Cref{thm: decomposition}, it effectively takes care of this kernel.

We see a mild dependency on $h$ for the larger values of $\alpha$. An explanation for this effect can be found in the proofs of \Cref{lem: stable transfer,lem: stable decomposition}. In \eqref{eq: transfer1 stable} and \eqref{eq: stability 2}, $\cb_k$ only appears in a product with $\alpha$. Thus, for small $\alpha$, the influence of the $h$-dependent constant is diminished.

The ASP-AMG solver outperforms the AMG solver, converging in all cases whereas the latter does not. We do observe a dependency on $h$, which may again be related to the Poincaré constant. Interestingly, there seems to be almost no dependency on $\alpha$. We acknowledge that AMG can benefit from a tailored smoother. However, such improvements can also be made for the steps in the ASP-AMG solver, and we leave such solver optimizations as future research. In conclusion, we note that the performance of an iterative solver improves if it is combined with the decomposition of \Cref{thm: decomposition}, particularly if the second-order differential terms dominate.

\section{Concluding remarks}
\label{sec: Concluding remarks}

In this work, we considered a particular decomposition of Hilbert complexes that leads to an explicit construction of a Poincaré operator. For our main example concerning the finite element de Rham complex, we showed that such a decomposition can be obtained by employing spanning trees in the grid. The construction handles the case of non-trivial cohomology classes which, for the case of the de Rham complex, means that the domain can have arbitrary topology.

The availability of a Poincaré operator led to several observations, including a second decomposition of the Hilbert complex. This decomposition forms a different basis for the Hilbert spaces, in which the mixed formulation of the Hodge-Laplace problem becomes a series of (at most) seven symmetric positive definite systems. By solving these systems in sequence, the solution to the original problem is recovered.
Numerically, this direct solver provides a significant speed-up compared to solving the saddle-point system in the original basis, without loss of accuracy.

The Poincaré operator moreover allows us to explicitly construct a basis for the spaces of harmonic forms. The numerically obtained harmonic forms are consistent with theoretical expectation.

Finally, we used the decomposition to propose auxiliary space preconditioners for mixed finite element problems. The numerical results showed that the preconditioner is robust for problems in which the differential operators dominate, by properly handling the large kernel of those terms.

The ideas developed in this work can directly be applied to other Hilbert complexes, including the discrete de Rham complex given by the virtual elements (\Cref{ex: VEM}). The identification of permitting decompositions for complexes such as the elasticity complex \cite{arnold2006finite} and the \v{C}ech-de Rham complex \cite{boon2025hodge}, forms a topic for future research.


\section*{Acknowledgments}

The author warmly thanks Alessio Fumagalli for the fruitful discussions preceding this work.

\section*{Data Availability Statement}

The research code associated with this article is available at \url{github.com/wmboon/poincare}.

\bibliographystyle{siam}
\bibliography{references}

\end{document}